\documentclass{interact}
\usepackage[utf8]{inputenc}

\usepackage{epstopdf}
\usepackage[caption=false]{subfig}
\usepackage{amsmath,amsfonts,amssymb,amsthm}
\usepackage{nicefrac}
\usepackage[numbers,sort&compress]{natbib}%
\bibpunct[, ]{[}{]}{,}{n}{,}{,}%
\makeatletter%
\def\NAT@def@citea{\def\@citea{\NAT@separator}}%
\makeatother%
\usepackage{multirow}

\usepackage{microtype}
\theoremstyle{plain}%
\newtheorem{theorem}{Theorem}[section]

\newtheorem{corollary}[theorem]{Corollary}
\newtheorem{proposition}[theorem]{Proposition}

\theoremstyle{definition}
\newtheorem{definition}[theorem]{Definition}
\newtheorem{example}[theorem]{Example}

\theoremstyle{remark}
\newtheorem{remark}{Remark}

\usepackage{subfig}
\usepackage{pgfplots,tikz}
\usetikzlibrary{shapes.geometric}
\pgfplotsset{compat=1.18}

\usepackage{booktabs}
\usepackage{multirow}

\usepackage[]{algorithm2e}
\usepackage[frozencache,cachedir=.]{minted}
\definecolor{bg}{rgb}{0.95,0.95,0.95}
\newlength\Colsep
\usepackage[colorlinks = true,
            linkcolor = blue,
            urlcolor  = blue,
            citecolor = teal,
            anchorcolor = purple]{hyperref}

\newcommand*\pFqskip{8mu}
\newcommand{\pFq}[5]{\mskip\pFqskip\relax{}_{#1}F_{#2}\biggl[ \begin{array}{c}#3\\#4\end{array} ;#5\biggr]}

\definecolor{mcyan}{rgb}{0.00000,0.44700,0.74100}

\title{Theoretical error estimates for computing the matrix logarithm by Padé-type approximants}
\author{\name{Lidia Aceto\textsuperscript{a} and Fabio Durastante\textsuperscript{b}\thanks{CONTACT F.~Durastante. Email: fabio.durastante@unipi.it}}
\affil{\textsuperscript{a} Dipartimento di Scienze e Innovazione Tecnologica, Università del Piemonte Orientale, Viale T. Michel, 11 - 15121 Alessandria, Italy; \textsuperscript{b} Dipartimento di Matematica, Università di Pisa, Via F. Buonarroti, 1/C - 56127 Pisa, Italy}}

\date{\today}
\begin{document}

\maketitle

\begin{abstract}
In this article, we focus on the error that is committed when computing the matrix logarithm using the Gauss--Legendre quadrature rules. These formulas can be interpreted as Padé approximants of a suitable Gauss hypergeometric function. Empirical observation tells us that the convergence of these quadratures becomes slow when the matrix is not close to the identity matrix, thus suggesting the usage of an inverse scaling and squaring approach for obtaining a matrix with this property. The novelty of this work is the introduction of error estimates that can be used to select {\it a priori} both the number of Legendre points needed to obtain a given accuracy and the number of inverse scaling and squaring to be performed. We include some numerical experiments to show the reliability of the estimates introduced.
\end{abstract}

\begin{keywords}
Matrix logarithm; Gauss--Legendre rule; Padé approximant; error bound.
\end{keywords}

\begin{amscode}
65F60; 47A58; 65D32; 41A20.
\end{amscode}

\section{Introduction}
In this work, we consider a problem that has a long history in Numerical Analysis, namely the computation of the matrix logarithm. A logarithm of 
a given non-singular matrix $A \in \mathbb{R}^{n \times n}$ is defined as any solution of the matrix equation \[e^X = A,\]
where $e^X$ denotes the matrix exponential of $X.$
When $A$ has no eigenvalues on the closed negative real axis, then there exists a unique real logarithm of $A$ whose eigenvalues lie in the strip $\{z \in  \mathbb{C}: -\pi < \mbox{Im}(z) < \pi \};$ see, e.g., \cite[Theorem 1.31]{Higham}. This unique logarithm is known as the {\it principal logarithm} of~$A$. Although we will later omit the word ``principal'' for simplicity, in this article we limit ourselves to this case as this is mostly used in a wide variety of applications and denote it with $\log(A)$.

Several approaches have been proposed in the literature for the evaluation of the matrix logarithm, but the most reliable and efficient seems to be the one that combines the Padé approximant method with an inverse scaling and squaring procedure~\cite{DieciMoriniPapini, KenneyLaub89}. This method first exploits the fact that 
\[
\log(A)= 2^s \log(A^{1/2^s})
\]
to determine a positive integer $s$ such that $A^{1/2^s}$ is close to the identity matrix and then computes $\log(A^{1/2^s})$ with a Padé approximant. Available algorithms using this strategy are~\cite{MR1825853,MR2970418,MR2200942} to be coupled with stable ways of evaluating the Padé approximant~\cite{MR1824061,MR3080997}. More recent results include variants using \emph{mixed-precision}~\cite{MR3775134}; the use of the dual inverse scaling and squaring
algorithm~\cite{MR4454938}, i.e., the solution of the matrix equation $r(X) = A^{2^{-s}}$ for $r(\cdot)$ a rational approximant to $e^z$ at $z = 0$; and the use of quadrature rules based upon a double exponential transformation~\cite{MR4074499}. For cases in which working with complex arithmetic is feasible, it is possible to use also quadrature formulas for contour integral formulations, see~\cite{MR2421045}.

For inverse scaling and squaring algorithms there are two points that should be addressed more systematically:
\begin{enumerate}
    \item the value of $s$ is typically obtained by progressively taking the square roots of $A$ up to $ ||A^{1/2^s}|| <1;$
    \item it is not possible to decide in advance which degree to use in the Padé approximant to obtain a given accuracy.
\end{enumerate}
The main goal of this work is therefore to provide a theoretical error analysis that leads to establishing these  parameters \emph{a priori}. The starting point for doing this is the following identity (cf. \cite[Eq. 15.4.1]{NIST:DLMF})
\begin{equation} \label{logGauss}
 \frac{\log(1+z)}{z} = \pFq{2}{1}{1,1}{2}{-z}
\end{equation}
which expresses the logarithmic function as a special case of the Gauss hypergeometric function. In fact, by virtue of this, the results on the Padé approximations to the Gauss hypergeometric functions help to find the results corresponding to $\log(1+z)/z,$  cf. \cite[Sec. 4]{Gomilko}.

The organization of the paper is as follows. In Section~\ref{sec2} 
we report some properties of hypergeometric functions useful for our analysis. In Section~\ref{sec3}, we prove that the $k$-point Gauss--Legendre rule used to approximate the function~\eqref{logGauss} is a $[k-1/k]$-Padé approximant. Starting from this fact, we present an explicit error estimate. The analysis is then extended in Section~\ref{sec4} to the case of the matrix logarithm. Some numerical experiments are proposed in Section~\ref{sec5}. Finally, the concluding remarks are given in Section~\ref{sec6}.

To have reproducible research, the code for generating all the results and figures is made available in the GitHub repository  \href{https://github.com/Cirdans-Home/padelogarithm}{Cirdans-Home/padelogarithm}.

\section{Padé approximants to Gauss hypergeometric functions} \label{sec2}

In this section, for the convenience of the reader, we recall the basic results for Gauss hypergeometric functions and report some properties of these functions useful for our analysis. The Gauss  hypergeometric function is meant the power series
\begin{equation} \label{gaussiper}
    \pFq{2}{1}{a,b}{c}{z}= \sum_{j=0}^{+\infty} \frac{(a)_j (b)_j}{(c)_j} \frac{z^j}{j!}, \qquad |z|<1,
\end{equation}
where $z$ is the complex variable,  $a, b, c$ are parameters that can take arbitrary real or complex values (provided that $c \neq 0, -1,-2,\dots$), and the symbol  $(q)_j$ denotes the Pochhammer symbol defined as
\[
{\displaystyle 
(q)_{j}=
\begin{cases}
1 & j=0,\\ 
q(q+1)\cdots (q+j-1) & \mbox{otherwise}.
\end{cases}
}
\]
When $ \mathfrak{Re}(c)> \mathfrak{Re}(b)> 0 $ we have an integral representation %
which constitutes an analytic continuation for~\eqref{gaussiper},   cf. \cite[Eq.~(15.3.1)]{Abra},
\begin{equation} \label{integrale}
    \pFq{2}{1}{a,b}{c}{z}= \frac{\Gamma(c)}{\Gamma(b) \Gamma(c-b)}\int_0^1 t^{b-1} (1-t)^{c-b-1} (1-tz)^{-a} \,\mathrm{d}t;
\end{equation}
here $\Gamma$ is the Euler gamma function, i.e., the analytic continuation to all complex numbers (except the non-positive integers) of the convergent improper integral function,
\begin{equation*}
    \Gamma(t) = \int_{0}^{+\infty}x^{t-1}e^{-x}\,\mathrm{d}x.
\end{equation*}
By setting in~\eqref{gaussiper} $b=1,$ the corresponding Gauss hypergeometric function can be approximated with the $[m/k]$-Padé approximant determined by the pair of polynomials $(P_{mk}^{[a,c]}(z),Q_{mk}^{[a,c]}(z))$.

Assuming that these two polynomials are normalized, that is  $Q_{mk}^{[a,c]}(0)=1,$  and that they are relatively prime,  as reported in \cite[eqs. (30)-(31)]{Gomilko}, for $m \ge k-1$ we have
\[
\pFq{2}{1}{a,1}{c}{z} \approx \frac{P_{mk}^{[a,c]}(z)}{Q_{mk}^{[a,c]}(z)},
\]
with
\begin{eqnarray}
    P_{mk}^{[a,c]}(z) &=& \sum_{j=0}^m \left( \sum_{\ell=0}^j 
    \frac{ (a)_{j-\ell} (-k)_\ell (-a-m)_\ell}{\ell! (1-c-m-k)_\ell (c)_{j-\ell}} \right) z^j, \label{num} \\
     Q_{mk}^{[a,c]}(z) &=& \pFq{2}{1}{-k,-a-m}{1-c-m-k}{z}. \label{den}
\end{eqnarray} 
Actually, to be more precise, the formula \eqref{num} for the numerator has been presented in \cite{Driver}, while Baker and Graves-Morris in \cite[p. 60]{BakG} give the formula \eqref{den} for the denominator (see also \cite[p. 65]{Bak}).

\section{The Gauss--Legendre approximation} \label{sec3}
As already said in the introduction, we start from the identity given in~\eqref{logGauss}, i.e.,
\[
\frac{\log(1+z)}{z} = \pFq{2}{1}{1,1}{2}{-z}.
\]
Thus, using~\eqref{integrale} and $\Gamma(1)=1$, we obtain an integral representation for this function which reads as
\[
\frac{\log(1+z)}{z} = \int_0^1  (1+tz)^{-1} \,\mathrm{d}t.
\]
Applying the change of variable $t \rightarrow (1+x)/2$ we can write
\begin{equation} \label{integleg}
    \frac{\log(1+z)}{z} = \int_{-1}^{1} \frac{1}{z(1+x)+2}\,\mathrm{d}x, \qquad z \in \mathbb{C} \setminus (-\infty,-1].
\end{equation}
The integrand is analytic in a neighborhood of $[-1,1]$ as long as $z$ stays away from $(-\infty,-1].$ Then, the Gauss--Legendre quadrature is an obvious choice for approximating this integral.  Denoting by $x_i$ and $\omega_i,$ $i=1,2,\dots,k,$ the nodes and weights of the $k$-point Gauss--Legendre quadrature on $[-1,1],$  we have that this quadrature is stated~as
\begin{equation}\label{eq:knodes}
    R_{k-1, k}(-z)= \sum_{i=1}^k \frac{\omega_i}{z(1+x_i)+2}.
\end{equation}
The most reliable and fastest way of computing the $\{x_i,\omega_i\}_{i=1}^{k}$ depends on the value of $k$, see, e.g., \cite[Fig. 2.1]{MR3033086}. For our purposes we generally need a low number of quadrature nodes, thus we rely on the default choice made by Chebfun's \mintinline{matlab}{legpts} routine~\cite{MR3033086}.

\begin{remark}
A completely different approach for attaining an integral expression of $\log(1+z)$ can be obtained by applying the transformation $x=-1-2/t$ to \eqref{integleg}, thus rewriting it as 
\begin{equation*}
    \frac{\log(1+z)}{z} = \int_{-\infty}^{-1} \frac{t^{-1}}{t-z}\,\mathrm{d}t.
\end{equation*}
This shows indeed that the function
\[
f(z) = \frac{\log(1+z)}{z},
\]
belongs to the class of Cauchy-Stieltjes functions. This allows obtaining an error analysis for the calculation of the matrix function based on the solution of the fourth problem of Zolotarev for which we refer to \cite{MR4235307}, and which is of use every time we know how to solve the Zolotarev problem for a set that contains the spectrum and generate the corresponding extremal rational function or its poles.

Let us also remark that for matrix Lie groups, such as the groups of orthogonal and symplectic matrices, there exist variants of the Gauss--Legendre approach discussed here that are shown to be structure-preserving, see~\cite{MR1826647}.
\end{remark}

\subsection{Not a diagonal Padé approximant}

Using the results reported at the end of Section~\ref{sec2} we know that the $[k-1/k]$-Padé approximant for the function given in \eqref{logGauss} reduces to
\begin{eqnarray} \label{padenostro}
     \pFq{2}{1}{1,1}{2}{-z} \approx  \frac{ P_{k-1,k}^{[1,2]}(-z)}{\displaystyle{\pFq{2}{1}{-k,-k}{-2k}{-z}}},
\end{eqnarray}
where 
\begin{equation} \label{numPD}
    P_{k-1, k}^{[1,2]}(-z) = \sum_{j=0}^{k-1} \left( \sum_{\ell=0}^{j} 
    \frac{ (1)_{j-\ell} (-k)_\ell (-k)_\ell}{\ell! (-2k)_\ell (2)_{j-\ell}} \right) (-z)^j.
\end{equation}
Now we want to determine what relationship exists between this Padé approximant  and the quadrature formula~\eqref{eq:knodes}.  If we can prove that the denominators of these two approximations coincide, we also know that the numerator is the same.  In fact, once the denominator in~\eqref{padenostro} is given, the numerator $P_{k-1,k}^{[1,2]}(-z)$ is determined by the definition of Padé approximant. 
Concerning the denominator, from \eqref{gaussiper} we have
\begin{eqnarray} \label{eq:1}
\pFq{2}{1}{-k,-k}{-2 k}{-z}&=&\sum_{j=0}^k  \frac{(-k)_j (-k)_j }{(-2k)_j} \frac{(-z)^j}{j!} \nonumber\\
&=& \sum_{j=0}^k \binom{k}{j}  \frac{(-k)_j }{(-2k)_j} z^j \nonumber \\
&=& \sum_{j=0}^k \binom{k}{j}  \frac{(k-j+1)_j }{(2k-j+1)_j} z^j \nonumber\\
&=& \frac{1}{\binom{2k}{k}}\sum_{j=0}^k \binom{k}{j} \binom{2k-j}{k}  z^j. 
\end{eqnarray}
Now, an explicit expression for the shifted Legendre polynomials is given by
\[
{\displaystyle {{\widetilde{\mathcal{L}}}_{k}(w)=(-1)^{k}\sum _{s=0}^{k}{\binom {k}{s}}{\binom {k+s}{s}}(-w)^{s}\,.}}
\]
Clearly, this implies that
\begin{eqnarray} \label{eq:2}
-z^k \widetilde{\mathcal{L}}_{k}(-z^{-1})&=&z^{k}\sum_{s=0}^{k}{\binom {k}{s}}{\binom {k+s}{s}}(z^{-1})^{s} \nonumber \\
&=& \sum_{s=0}^{k}{\binom {k}{k-s}}{\binom {k+s}{k}}z^{k-s} \nonumber \\
 &=& \sum_{j=0}^{k}{\binom {k}{j}}{\binom {2k-j}{k}}z^j.
\end{eqnarray}
Therefore, from~\eqref{eq:1} and~\eqref{eq:2} we deduce the relation that connects the denominator of the Padé approximant with the shifted Legendre polynomials, namely
\[
   \pFq{2}{1}{-k,-k}{-2 k}{-z} =  \frac{1}{\binom{2k}{k}} (-z)^k {\widetilde{\mathcal{L}}}_{k}(-z^{-1}).
\]
However, considering that the shifted Legendre polynomials are defined in terms of Legendre polynomials as
\[
{\displaystyle {\widetilde{\mathcal{L}}}_{k}(w)={\mathcal{L}}_{k}(2w-1)\,,}
\]
it is also immediate to obtain the following identity which relates the denominator of the approximant Padé with the Legendre polynomials
\begin{equation} \label{goodden}
     \pFq{2}{1}{-k,-k}{-2 k}{-z} =  \frac{1}{\binom{2k}{k}} (-z)^k {\mathcal{L}}_{k}\left(-\frac{2}{z}-1\right).
\end{equation}
It is worth pointing out that this relation also expresses the connection we were looking for. In fact, multiplying  by $\binom{2k}{k}$ we find that the right side thus obtained coincides with the denominator of the rational function which defines the Gauss--Legendre quadrature formula in~\eqref{eq:knodes}. Therefore, the following result holds.
\begin{proposition}\label{pro:padeformulation}
  The Gauss--Legendre rule in~\eqref{eq:knodes} is the $[k-1/k]$-Padé approximant for the function $\log(1+z)/z$ given by
  \begin{equation*}
      R_{k-1, k}(-z)=  \frac{\binom{2k}{k} P_{k-1,k}^{[1,2]}(-z)}{\binom{2k}{k}\pFq{2}{1}{-k,-k}{-2 k}{-z}},
  \end{equation*}
  where $P_{k-1,k}^{[1,2]}(-z)$ is defined by~\eqref{numPD}.
\end{proposition}

\begin{remark}
We stress that the connection between Gauss quadratures and Padé approximants is a well-known result, see, e.g., \cite[p.~34]{MR561106}. For the case of the Gauss--Legendre formula, the observation is available in \cite[Lemma 4.6]{FrommerGuttel} without deriving the expression given by Proposition~\ref{pro:padeformulation}. {We observe also that $z R_{k-1, k}(-z)$ is then the $[k,k]$-Padé approximant of the $\log(1+z)$ function, see, e.g., \cite[Section~11.4]{Higham}.}
\end{remark}

\subsection{Error analysis}
In order to obtain an estimate of the error for the rational approximation $z R_{k-1,k}(-z)$, we consider results based on the theory of Gauss hypergeometric functions. In fact, following \cite[Theorem 2]{Driver} we know that
\begin{equation*}
  \pFq{2}{1}{1,1}{2}{-z} - \frac{P_{k-1,k}^{[1,2]}(-z)}{\pFq{2}{1}{-k,-k}{-2 k}{-z}} =  \frac{\pFq{2}{1}{k+1,k+1}{2k+2}{-z}}{\pFq{2}{1}{-k,-k}{-2 k}{-z}} \frac{z^{2k}}{2k+1} \binom{2k}{k}^{-2}.
\end{equation*} 
Then, from the above considerations, we obtain 
\[
\log(1+z) - z \, R_{k-1,k}(-z) =  
\frac{\pFq{2}{1}{k+1,k+1}{2k+2}{-z}}{\pFq{2}{1}{-k,-k}{-2 k}{-z}} \frac{z^{2k+1}}{2k+1} \binom{2k}{k}^{-2}.
\]
Using \eqref{goodden} and recalling that the Legendre polynomials have definite parity, that is they are even or odd according to  
\[
{\displaystyle {{\mathcal{L}}}_{k}(-w)=(-1)^{k} {{\mathcal{L}}}_{k}(w),}
\]
the previous relationship becomes
\begin{equation} \label{toapprox}
\log(1+z) - z \, R_{k-1,k}(-z) = 
  \frac{\pFq{2}{1}{k+1,k+1}{2k+2}{-z}}{{\mathcal{L}}_{k}\left(\frac{2}{z}+1 \right)} \frac{z^{k+1}}{2k+1}\binom{2k}{k}^{-1}.
\end{equation}
In the form in which it is found, this result is not very convenient for numerical work. However, if we assume that $k$ is large, then we may obtain an asymptotic estimate for the error in terms of elementary functions. To this aim, we use the  asymptotic form of the Gauss hypergeometric functions given in \cite[eq. (16) p. 77]{Erdelyi}.
\begin{eqnarray*}
  &&  (z/2-1/2)^{-a-\lambda} \pFq{2}{1}{a+\lambda, a-c+1+\lambda}{a-b+1+2\lambda}{2(1-z)^{-1}} = \\ 
  &=&\frac{2^{a+b} \Gamma(a-b+1+2\lambda) \Gamma(1/2) \lambda^{-1/2}}{\Gamma(a-c+1+\lambda) \Gamma(c-b+\lambda)} e^{-(a+\lambda) \xi} \\
    &\times& (1-e^{-\xi})^{-c+1/2} (1+e^{-\xi})^{c-a-b-1/2} [1+ \mathcal{O}(\lambda^{-1})]
\end{eqnarray*}
Such a formula applied to the numerator in~\eqref{toapprox} gives:
\begin{eqnarray*}
  &&  \left(\frac{v-1}{2}\right)^{-1-k} \pFq{2}{1}{k+1, k+1}{2k+2}{\frac{2}{1-v}} = \frac{2 \Gamma(2k+2) \Gamma(1/2) k^{-1/2}}{\Gamma(k+1) \Gamma(k+1)} \\
    &\times& e^{-(k+1) \xi}  (1-e^{-\xi})^{-1/2}(1+e^{-\xi})^{-1/2}    [1+ \mathcal{O}(k^{-1})]\\
    &=& 2 (2k+1) \binom{2k}{k}\sqrt{\frac{\pi}{k}} \left(v - (v^2-1)^{1/2}\right)^{k+1} \\
    &\times&   2^{-1/2}\left(1-v \left(v - (v^2-1)^{1/2}\right)\right)^{-1/2}   [1+ \mathcal{O}(k^{-1})],
\end{eqnarray*}
where 
\[
v - (v^2-1)^{1/2} = e^{- \xi}.
\]
Considering that 
$
{\displaystyle \Gamma \left({ {1}/{2}}\right)={\sqrt {\pi }}}
$
and setting
\[
\frac{2}{1-v} = -z
\]
we obtain
\begin{eqnarray}
  &&  z^{k+1} \pFq{2}{1}{k+1, k+1}{2k+2}{-z} =  (2k+1) \binom{2k}{k}\sqrt{\frac{2\pi}{k}} \left( 1 +\frac{2}{z}  -\frac{2}{z}(1+z)^{1/2} \right)^{k+1}  \nonumber \\
    &\times& 
     \left(1 -\left(1+\frac{2}{z}\right) \left( 1 +\frac{2}{z}  -\frac{2}{z}(1+z)^{1/2} \right)\right)^{-1/2}
    [1+ \mathcal{O}(k^{-1})]. \label{eq:ratio-numerator}
\end{eqnarray}
Instead, for the denominator in~\eqref{toapprox} we observe that 
\[
{{\mathcal{L}}}_k(w) \sim \frac{1}{\sqrt{2 \pi k}} \frac{(w+(w^2-1)^{1/2})^{k+1/2}}{(w^2-1)^{1/4}}.
\]
Therefore, setting $w = \nicefrac{2}{z} +1$, we find
\begin{equation}
{{\mathcal{L}}}_k \left(\frac{2}{z}+1 \right) \sim \frac{1}{\sqrt{2 \pi k}} \frac{(\frac{2}{z}+1+ \frac{2}{z}(1+z)^{1/2})^{k+1/2}}{\left(\frac{2}{z}( 1+z)^{1/2}\right)^{1/2}}.\label{eq:ratio-denominator}
\end{equation}
The ratio of the two asymptotic formulas for the numerator~\eqref{eq:ratio-numerator} and the denominator~\eqref{eq:ratio-denominator} -- having multiplied both of them by $\left( 1 +\frac{2}{z}  +\frac{2}{z}(1+z)^{1/2} \right)^{1/2}$ -- lead to the following relation which gives the expression for the remainder  
\begingroup
\allowdisplaybreaks
\begin{eqnarray} \label{ratio}
\log(1+z) - z \, R_{k-1,k}(-z) 
&\sim&  2\pi \left( 1 +\frac{2}{z}  +\frac{2}{z}(1+z)^{1/2} \right)^{1/2}    \left(\frac{2}{z} (1+z)^{1/2} \right)^{1/2}\nonumber\\
    &\times&\left(1 -\left(1+\frac{2}{z}\right) \left( 1 +\frac{2}{z}  -\frac{2}{z}(1+z)^{1/2} \right)\right)^{-1/2}  
    \nonumber\\
   &\times&  \left( \frac{ 1+\frac{2}{z} - \frac{2}{z}(1+z)^{1/2} }{ 1+\frac{2}{z} + \frac{2}{z}(1+z)^{1/2} }\right)^{k+1}  \nonumber\\  
 &=&     2\pi \left( \frac{2}{z} (1+z)^{1/2}  +\frac{4}{z^2}(1+z)^{1/2} +\frac{4}{z^2}(1+z) \right)^{1/2}    \nonumber\\
    &\times&\left( \frac{2}{z}(1+z)^{1/2}
+   \frac{4}{z^2} (1+z)^{1/2} -\frac{4}{z^2}(1+z) \right)^{-1/2}  
    \nonumber\\
   &\times&  \left( \frac{ 1+\frac{2}{z} - \frac{2}{z}(1+z)^{1/2} }{ 1+\frac{2}{z} + \frac{2}{z}(1+z)^{1/2} }\right)^{k+1}.
 \end{eqnarray}
\endgroup
Consequently, \eqref{ratio} becomes
\begingroup
\allowdisplaybreaks
\begin{eqnarray} \label{ratioLAST}
&&\log(1+z) - z \, R_{k-1,k}(-z) 
\sim  2\pi  
 \left( \frac{ 1+\frac{2}{z} + \frac{2}{z}(1+z)^{1/2} }{ 1+\frac{2}{z} - \frac{2}{z}(1+z)^{1/2} }\right)^{1/2} \nonumber\\
&\times&
\left( \frac{ 1+\frac{2}{z} - \frac{2}{z}(1+z)^{1/2} }{ 1+\frac{2}{z} + \frac{2}{z}(1+z)^{1/2} }\right)^{k+1}  = 2\pi  \left( \frac{ 1+\frac{2}{z} - \frac{2}{z}(1+z)^{1/2} }{ 1+\frac{2}{z} + \frac{2}{z}(1+z)^{1/2} }\right)^{k+\frac{1}{2}}\nonumber \\
&=& 2\pi  \left( \frac{ z+ 2  - 2(1+z)^{1/2} }{ z+2 + 2(1+z)^{1/2} }\right)^{k+\frac{1}{2}} =2\pi  \left( \frac{1+ (1+z)  - 2(1+z)^{1/2} }{ 1+ (1+z) + 2(1+z)^{1/2} }\right)^{k+\frac{1}{2}}\nonumber \\
&=& 2 \pi \left( \frac{(1- (1+z)^{1/2})^2   }{ (1+ (1+z)^{1/2})^2 }\right)^{k+\frac{1}{2}} \nonumber\\
&=&2 \pi \left( \frac{1- (1+z)^{1/2}   }{ 1+ (1+z)^{1/2}}\right)^{2k+1}.
\end{eqnarray}
\endgroup
This concludes the analysis, which we can summarize in the following result.
\begin{proposition}\label{pro:scalar-bound}
{The error with respect to the $[k/k]$-Padé approximation  of $\log(1+z)$ based on the approximation in~\eqref{eq:knodes}~is}
\[
|\log(1+z) - z \, R_{k-1,k}(-z)| \sim 2 \pi \left\lvert \left( \frac{1- (1+z)^{1/2}   }{ 1+ (1+z)^{1/2}}\right)^{2k+1} \right\rvert, \quad k \rightarrow + \infty.
\]
\end{proposition}
\begin{example}\label{example:goodnessofthescalarbound}
We can illustrate the result of Proposition~\ref{pro:scalar-bound} by comparing the error with the bound for increasing values of $k$ with a few lines of \texttt{MATLAB} code
\begin{minted}[fontsize=\small]{matlab}
k = 3; %
[x,omega] = legpts(k);
R = @(z,x,omega) arrayfun(@(t) sum(omega./(t.*(1+x) + 2)),-z);
bound = @(z,k) arrayfun(@(t) max(2*pi*abs(((1 - sqrt(1+t))...
    ./(1+sqrt(1+t))).^(2*k+1)),eps),z);
err = @(z) abs(log(1+z) - z.*R(-z,x,omega.'));
z = linspace(-1,3,100);
semilogy(z,err(z),'-',z,bound(z,k),'--','LineWidth',2)
grid on, axis tight, title(sprintf('k = %
legend({'Error','Estimate'})
\end{minted}
The repeated application for $k=3,6,\ldots,18$ returns the results in Figure~\ref{fig:ratio_bound} which show a very good agreement between the estimate and the actual error.

\begin{figure}[htbp]
    \centering
    \input{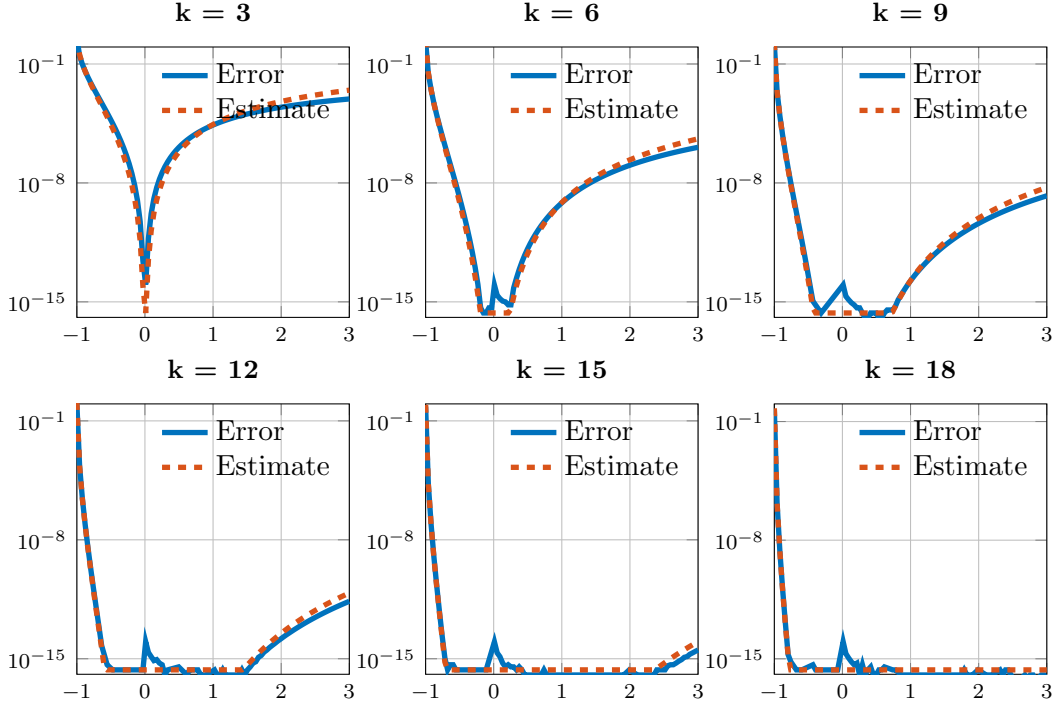}
    \caption{Example \ref{example:goodnessofthescalarbound}. Depiction of the bound from Proposition~\ref{pro:scalar-bound} for the Gauss--Legendre approximation for different values of $k$ in the $(-1,3]$ interval.}
    \label{fig:ratio_bound}
\end{figure}

\end{example}

\section{Using the estimate for the matrix logarithm} \label{sec4}

The bound from Proposition~\ref{pro:scalar-bound} is a bound for scalar functions thus we need to transform it into a bound for matrix functions, i.e., we need a way to estimate
\[
\| \log(I + B) - B R_{k-1,k}(-B) \| = \| E_k(B) \|,
\]
for a given matrix $B \in \mathbb{R}^{n \times n}$ such that $I + B$ has no eigenvalues on the closed negative real axis and being $\|\cdot\|$ the $2$-norm. For matrices $B$ such that $\|B\| < 1$ the common procedure to estimate the bound is using that for any Padé $R_{m,k}$ approximant~\cite{KenneyLaub89} 
\[
\| \log(I + B) - R_{m,k}(-B)\| \leq \lvert \log(I + \|B\|) - R_{m,k}(-\|B\|) \rvert.
\]
Indeed, under this restriction, and selecting the rational approximant as in~\eqref{eq:knodes}, we can immediately obtain a bound using Proposition~\ref{pro:scalar-bound}.
\begin{corollary}\label{cor:bound-with-restriction}
Let $B$ be a matrix for which $\log(I + B)$ admits a principal logarithm with $\|B\| < 1$, then 
\[
\| \log(I + B) - B R_{k-1,k}(-B)\| \sim 2 \pi \left\lvert \left( \frac{1- (1+\|B\|)^{1/2}   }{ 1+ (1+\|B\|)^{1/2}}\right)^{2k+1} \right\rvert, \quad k \rightarrow + \infty.
\]
\end{corollary}

\begin{figure}[htbp]
    \centering
    \includegraphics[width=\columnwidth]{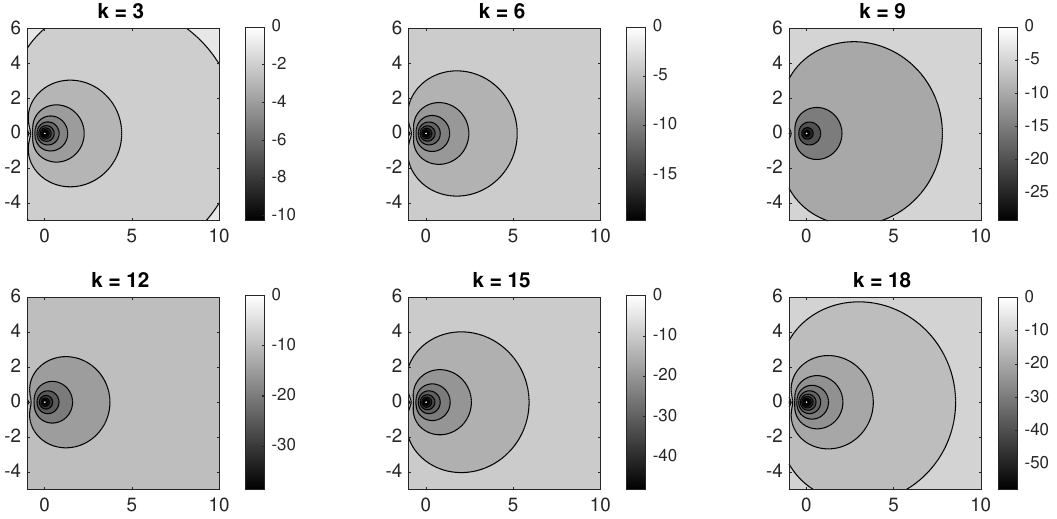}
    \caption{Level sets in $\log_{10}(\cdot)$ scale of the absolute value of~\eqref{ratioLAST} for different values of $k$.}
    \label{fig:complex-bound}
\end{figure}

To remove such a constraint on the norm of $B$ we need to use the concept of \emph{field of values} of a matrix.
\begin{definition}
  Given an $A \in \mathbb{C}^{n \times n}$ we define its \emph{field
    of values} as the subset of the complex plane 
  \[
    \mathcal W(A) = \left\lbrace \left.\frac{\mathbf{x}^H A \mathbf{x}}{\mathbf{x}^H \mathbf{x}} \ \right| \ \mathbf{x} \in \mathbb {C}^n \setminus \{ 0 \} \right\rbrace. 
  \]
\end{definition}
For a normal matrix $A$, $\mathcal{W}(A)$ reduces to the convex hull of the eigenvalues, while in the general case one only knows that the spectrum is contained in $\mathcal W(A)$, with the latter possibly being significantly larger. This is of interest to us because it permits us to state the following results.
\begin{proposition}[{Crouzeix and Palencia~\cite{MR3666309}}]\label{pro:crouzeix}
  Let $A \in \mathbb{C}^{n \times n}$ be a matrix, and $g(z)$ be a holomorphic
  function defined on $\mathcal W(A)$. Then,
  \[
    \|g(A)\| \leq \mathcal{C} \max_{x \in \mathcal W(A)} |g(x)|, 
  \]
  where $\mathcal{C}$ is a universal constant smaller or equal than $1+\sqrt 2$. 
\end{proposition}
This let us discharge the problem of bound $E_k(B)$ to the problem of bounding the maximum of $E_k(x)$ on $\mathcal{W}(B)$ whenever $\mathcal{W}(B) \subset \{ x \in \mathbb{C} \,:\,\mbox{Re}(x) > -1 \}$, that is
\begin{equation}\label{eq:bound-on-the-complex}
\begin{split}
\| \log(I + B) - B R_{k-1,k}(-B) \| \leq &\; \mathcal{C} \max_{x \in \mathcal{W}(B)} |\log(1+x) - x R_{k-1,k}(-x)| \\ \leq &\;  2 (1+\sqrt{2})\pi \max_{x \in \mathcal{W}(B)}  \left\lvert \left( \frac{1- (1+x)^{1/2}   }{ 1+ (1+x)^{1/2}}\right)^{2k+1} \right\rvert.
\end{split}
\end{equation}
Due to the maximum modulus principle -- \eqref{ratioLAST} is holomorphic on $\mathcal{W}(B)$ -- then the maximum of its absolute value in~\eqref{eq:bound-on-the-complex} is reached on the boundary of $\mathcal{W}(B)$, see, e.g., the plots of the level sets in logarithmic scale in Figure~\ref{fig:complex-bound}.
This gives us a precise description of the behavior that has been observed many times in the literature, regarding the accuracy of using~\eqref{eq:knodes} with spectra with eigenvalues with an absolute value larger than $1$, while also expelling the deterioration due to the eigenvalues of $B$ with the real part near to $-1$; see again the left part of the estimate in Figure~\ref{fig:ratio_bound}.

\begin{example}\label{example:daisyexample}
To illustrate the bound based on the field of values, we consider a sequence of Toeplitz matrices, i.e., of matrices that are constant along the diagonals, given by
\[
(T_n)_{i,j} = \begin{cases}
2.5, & i-j = 0, \\
-1, & i-j=1,\\
1, & i-j=-5,
\end{cases} \quad i,j=0,\ldots,n-1,
\]
we associate to $T_n \in \mathbb{R}^{n \times n}$ its generating function $a(\theta) = 2.5 - e^{i\theta} + e^{-5i\theta}$, that is simply the function whose Fourier coefficients give the entries of $(T_n)_{i,j}$. Generating functions convey information about the infinite operator of whose the matrix $T_n$ is a finite section, and can be used to infer several spectral information regarding the matrix sequence, see, e.g, Figure~\ref{fig:daisyexample-panela}.
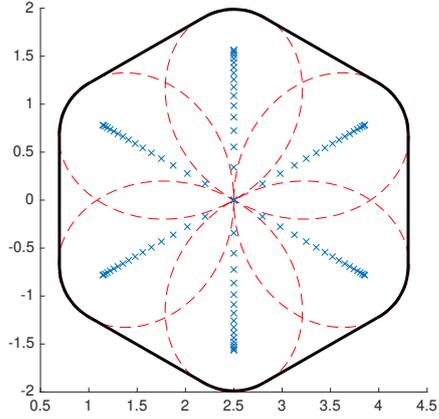
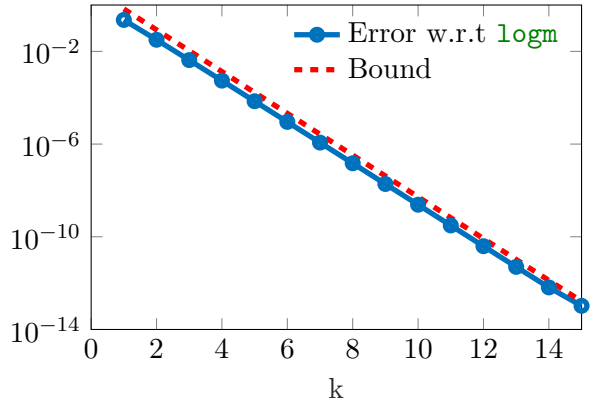
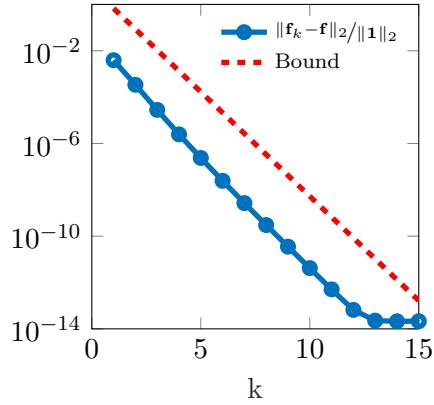
\begin{figure}[htbp]
    \centering
    \subfloat[Generating function (\raisebox{0.3em}{\protect\tikz{\protect\draw[-,dashed,red] (0,0)--(0.38,0);}}), eigenvalues ($\textcolor{mcyan}{\times}$) and field of values (\raisebox{0.3em}{\protect\tikz{\protect\draw[-,thick,black] (0,0)--(0.38,0);}}).\label{fig:daisyexample-panela}]{\includegraphics[width=0.4\columnwidth]{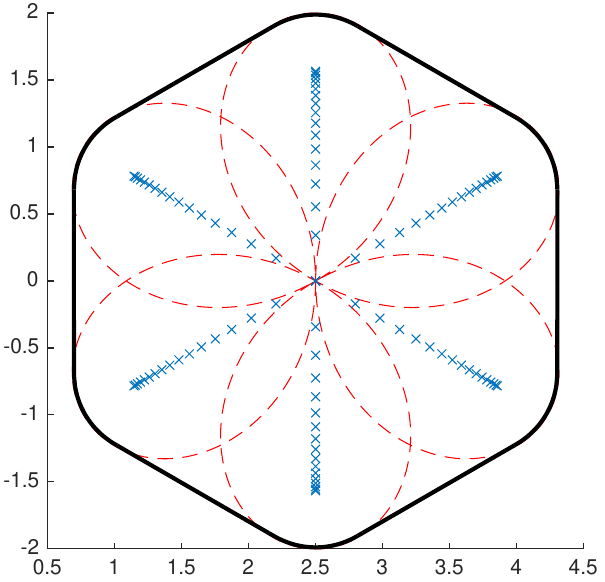}}\hfill
    \subfloat[Absolute error of the approximation based on~\eqref{eq:knodes} (\raisebox{0em}{\protect\tikz{\protect\draw[-,thick,mcyan] (0,0)--(0.38,0);\protect\filldraw[mcyan] (0.19,0) circle (0.3em);}}) and theoretical bound~\eqref{eq:bound-on-the-complex} (\raisebox{0.3em}{\protect\tikz{\protect\draw[-,dashed,thick,red] (0,0)--(0.38,0);}}) \label{fig:daisyexample-panelb} for $T_{500}$.]{% This file was created by matlab2tikz.
\definecolor{mycolor1}{rgb}{0.00000,0.44700,0.74100}%
\begin{tikzpicture}

\begin{axis}[%
width=0.45\columnwidth,
height=0.3\columnwidth,
at={(0in,0in)},
scale only axis,
xmin=0,
xmax=15,
xlabel style={font=\color{white!15!black}},
xlabel={k},
ymode=log,
ymin=1e-14,
ymax=1,
yminorticks=true,
axis background/.style={fill=white},
legend style={legend cell align=left, align=left, draw=none, fill=none}
]
\addplot [color=mycolor1, line width=2.0pt, mark=o, mark options={solid, mycolor1}]
  table[row sep=crcr]{%
1	0.226884195264898\\
2	0.0316657707908958\\
3	0.00419994850817343\\
4	0.000546169286632361\\
5	7.04263215091495e-05\\
6	9.04215313207797e-06\\
7	1.15803174742742e-06\\
8	1.48072933368406e-07\\
9	1.89130115374593e-08\\
10	2.41385866127209e-09\\
11	3.07909107477967e-10\\
12	3.92657253611775e-11\\
13	5.01264871713158e-12\\
14	6.52765301424244e-13\\
15	1.04903735951256e-13\\
};
\addlegendentry{Error w.r.t \mintinline{matlab}{logm}}

\addplot [color=red, dashed, line width=2.0pt]
  table[row sep=crcr]{%
1	0.675434186363572\\
2	0.0848531614521734\\
3	0.0106598972243211\\
4	0.00133917707824164\\
5	0.000168237573885431\\
6	2.11352790655732e-05\\
7	2.65517393566228e-06\\
8	3.33563072753738e-07\\
9	4.19047211975448e-08\\
10	5.26438866313111e-09\\
11	6.61352401459889e-10\\
12	8.30841009099463e-11\\
13	1.04376544317014e-11\\
14	1.31125725430547e-12\\
15	1.64730074004611e-13\\
};
\addlegendentry{Bound}

\end{axis}
\end{tikzpicture}%
}
    
    \subfloat[Bound for $\mathbf{f} = \log(T_{500})\mathbf{1}$ and $\mathbf{f}_k$ the approximation in the rational Krylov space with $k$ poles~\eqref{eq:knodes}.\label{fig:daisyexample-krylov}]{% This file was created by matlab2tikz.
\definecolor{mycolor1}{rgb}{0.00000,0.44700,0.74100}%
\begin{tikzpicture}

\begin{axis}[%
width=0.3\columnwidth,
height=0.3\columnwidth,
at={(0.772in,0.473in)},
scale only axis,
xmin=0,
xmax=15,
xlabel style={font=\color{white!15!black}},
xlabel={k},
ymode=log,
ymin=1e-14,
ymax=1,
yminorticks=true,
axis background/.style={fill=white},
legend style={legend cell align=left, align=left, draw=none, fill=none, font=\footnotesize}
]
\addplot [color=mycolor1, line width=2.0pt, mark=o, mark options={solid, mycolor1}]
  table[row sep=crcr]{%
1	0.00390777594414222\\
2	0.000339427000622663\\
3	2.80957671044206e-05\\
4	2.46237205527493e-06\\
5	2.36496650183965e-07\\
6	2.44316493212789e-08\\
7	2.66154407287756e-09\\
8	3.01350076121274e-10\\
9	3.50986921074475e-11\\
10	4.17436471473952e-12\\
11	5.05268171020071e-13\\
12	6.58162812688611e-14\\
13	2.26259592072669e-14\\
14	2.1201285991178e-14\\
15	2.14025175403291e-14\\
};
\addlegendentry{$\nicefrac{\|\mathbf{f}_k - \mathbf{f}\|_2}{\|\mathbf{1}\|_2}$}

\addplot [color=red, dashed, line width=2.0pt]
  table[row sep=crcr]{%
1	0.675434186363572\\
2	0.0848531614521734\\
3	0.0106598972243211\\
4	0.00133917707824164\\
5	0.000168237573885431\\
6	2.11352790655732e-05\\
7	2.65517393566228e-06\\
8	3.33563072753738e-07\\
9	4.19047211975448e-08\\
10	5.26438866313111e-09\\
11	6.61352401459889e-10\\
12	8.30841009099463e-11\\
13	1.04376544317014e-11\\
14	1.31125725430547e-12\\
15	1.64730074004611e-13\\
};
\addlegendentry{Bound}

\end{axis}
\end{tikzpicture}%
}
    
    \caption{Example~\ref{example:daisyexample}. In panel~\ref{fig:daisyexample-panela} we report the spectral information obtained through the generating function $a(\theta)$, in panel~\ref{fig:daisyexample-panelb} we report instead the absolute error evaluated with respect to \texttt{MATLAB}'s \mintinline{matlab}{logm} function against the error bound. Finally, panel~\ref{fig:daisyexample-krylov} gives a comparison of the bound for the rational Krylov method with poles given by~\eqref{eq:knodes}.}
\end{figure}
We compute $\log(T_{500})$ by means of~\eqref{eq:knodes}, and evaluate the bound by using as an approximation of the field of value the maximum and the minimum of $\lvert a(\theta)\rvert$ on $[-\pi,\pi]$. To have a reference logarithm we use \texttt{MATLAB}'s \mintinline{matlab}{logm} function. From the results depicted in Figure~\ref{fig:daisyexample-panelb}, we find a good accordance between the two quantities. We observe in passing that in this case, the bound from Corollary~\ref{cor:bound-with-restriction} is not applicable since $\|T_{500} - I\| \approx 3.4015 > 1$. Observe also that this permits us to obtain a bound for the computation of $\log(T_{500})\mathbf{v}$, $\mathbf{v} \in \mathbb{R}^{500}$, with a rational Krylov method~\cite{FrommerGuttel} since we can run it with the poles
\[
\xi_\ell = -\frac{2}{1 + x_\ell}, \qquad \ell=1,2,\ldots,k
\]
obtained from~\eqref{eq:knodes}; see, e.g., the result in Figure~\ref{fig:daisyexample-krylov} in which we compare
\[
    \| V_k \log(I_k + V_k^T (T_{500} - I) V_k) V_k^T \mathbf{1} - \log(T_{500})\mathbf{1}\|, \quad k=1,\ldots,15, \quad \mathbf{1} = (1,1,\ldots,1)^T,
\]
for $V_k$ the basis of the rational Krylov space $\mathcal{K}_k(T_{500}-I,\mathbf{v},\{\xi_\ell\}_{\ell=1}^{k})$, and the bound~\eqref{eq:bound-on-the-complex}.
\end{example}

A limitation of the proposed bound is given by the request to have a matrix whose field-of-value is contained in the right half plane. There are matrices that, despite having eigenvalues in the right half plane, have a field of values that contains the origin, e.g.,
\begin{equation}\label{eq:fovfailure}
A = \begin{bmatrix}
0.1 & 10^{6} \\
0 & 0.1
\end{bmatrix}, \qquad \mathcal{W}(A) = \vcenter{\hbox{\input{badfov.tikz}}}.
\end{equation}
An alternative in this case is represented by the usage of the $\epsilon$-pseudospectrum of the matrix $A$.
\begin{definition}
We define the $\epsilon$-pseudospectra as
\[
\Lambda_\epsilon(A) = \{ z \in \mathbb{C} \,:\, \|(zI - A)^{-1}\| \geq \epsilon^{-1}\}.
\]
\end{definition}
For which one can prove a result analogous to Proposition~\ref{pro:crouzeix}.
\begin{proposition}\label{pro:pseudospectra}
Let $f$ be a function that is analytic on $\Lambda_\epsilon(A)$ for a fixed $\epsilon > 0$. Then, provided that the boundary $\partial \Lambda_\epsilon(A)$ of $\Lambda_\epsilon(A)$ consists of a finite union of Jordan curves,
\[
\|f(A)\| \leq \frac{\mathcal{L}(\partial \Lambda_\epsilon(A))}{2\pi\varepsilon} \max_{z \in \Lambda_\epsilon(A)} |f(z)|,
\]
for $\mathcal{L}(\partial \Lambda_\epsilon(A))$ the arc-length of the boundary.
\end{proposition}

\begin{figure}[htbp]
    \centering
    \subfloat[Level set]{\includegraphics[width=0.45\columnwidth]{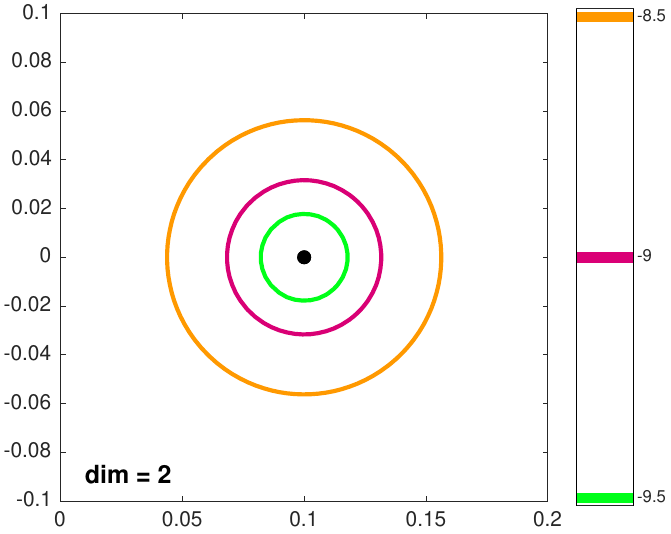}}\hfil
    \subfloat[Resolvent norm]{\includegraphics[width=0.45\columnwidth]{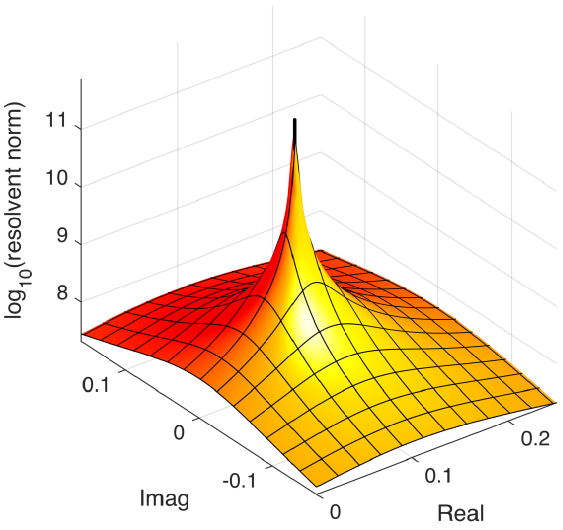}}
    \caption{$\epsilon$-pseudospectra of the matrix $A$ from~\eqref{eq:fovfailure} computed with \texttt{EigTool} (\href{https://github.com/eigtool/eigtool}{github.com/eigtool/eigtool}).}
    \label{fig:pseudospectra}
\end{figure}

For the example given in~\eqref{eq:fovfailure}, we find that the $\epsilon$-pseudospectra for $\epsilon = 10^{-8.5}$ and below are contained in the positive semi-plane, see Figure~\ref{fig:pseudospectra}. 

An algorithm for obtaining a numerical evaluation of the field of values is available in~\cite{MR474755} and implemented in Chebfun's \mintinline{matlab}{fov} command. Furthermore, the possibility and manner of constructing other spectral sets that can be used to obtain inequalities such as the one in Proposition~\ref{pro:crouzeix} has recently been investigated, see, e.g.,~\cite{MR4002721}. This could be of use to obtain spectral sets smaller than the field of values and which, above all, exclude the branch cut of the logarithm.
{On the other hand, the computation, i.e. the approximation, of the $\epsilon$-pseudospectra can be significantly more difficult than the computation of the field of value~\cite{MR1819647}. Typically there are two approaches, that of evaluating the resolving function on a grid of points and then considering the level sets as the edges of the regions of interest~\cite{MR1410089,MR1433795}--albeit with some tricks to reduce the computational intensity of the task--and that of using algorithms based on curve tracing~\cite{MR1410090}. Some algorithms have also been devised for the large-scale case~\cite{MR1861267,MR1889788}.}

Each of these bounds should be accompanied by an analysis of the accuracy that is actually achievable, for the implementation of the formula~\eqref{eq:knodes}. This is expected to be limited in norm by the product of the conditioning of the linear systems to be solved
\[
\begin{split}
\kappa = &\; \max_{i=1,\ldots,k} \kappa_2((2I + (1+x_i)B)/\omega_i) \\ = &\; \max_{i=1,\ldots,k} \| (2I + (1+x_i)B)/\omega_i \| \|\omega_i(2I + (1+x_i)B)^{-1}\|,
\end{split}
\]
and the unit round-off in the given precision\footnote{E.g., for \texttt{MATLAB} in double precision this is $u = 2^{-53} \approx 10^{-16}$.} times a constant depending only on the order of the rational approximation $k$; see~\cite[Section~3]{MR1824061}.

\subsection{Inverse scaling and squaring}

In general, the spectrum and the field of values of a matrix for which we want to compute the logarithm will not be in the region of rapid convergence described by the bound. The procedure in these cases is to transport the spectrum to a region where the Padé expansion is a better approximation using an \emph{inverse scaling and squaring algorithm}, i.e., exploit the property of the logarithm for which
\[
\log(A)= 2^s \log(A^{1/2^s}), \quad s \in \mathbb{N},
\]
and observing that under the same assumptions on $A$ for which a principal logarithm exists, there is a unique matrix $X$ satisfying $X^2 = A$, i.e., the principal square root of~$A$. The key point of this procedure is clearly the identification of the appropriate value of $s$ to reach an approximation of the logarithm within a predetermined tolerance.

\begin{figure}[htbp]
    \centering
    \input{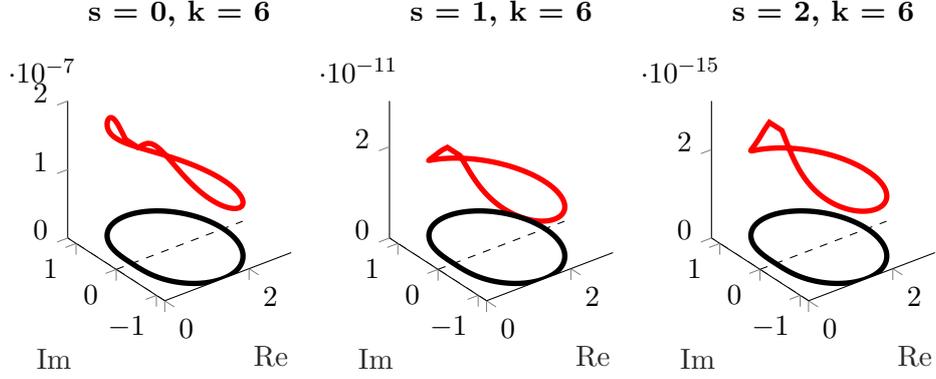}
    \caption{Depiction of the bound~\eqref{eq:findingthevalues} for the $10 \times 10$ matrix \mintinline{matlab}{A = gallery('forsythe',10,1e-10,0);}. The thick black line on the $x-y$ plane is the field of value of $A$, while the red line is the bound \eqref{eq:findingthevalues} for different values of $s$ and $k$.}
    \label{fig:3dbound}
\end{figure}

Since we compute
\[
\begin{split}
\log(A) = &\; 2^s \log(A^{1/2^s}) = 2^s \log(I + (A^{1/2^s} - I)) \\ = &\; 2^s (A^{1/2^s} - I)\sum_{i=1}^k \omega_i \left(2I + (1+x_i)(A^{1/2^s} - I)\right)^{-1} \\ = &\; 
2^s (A^{1/2^s} - I)\sum_{i=1}^k \omega_i \left((1 - x_i)I + (1+x_i)A^{1/2^s}\right)^{-1},
\end{split}
\]
we need to work with the bounds given by
\[
\|E(A,s,k)\| = \|E_k(A^{1/2^s} - I)\| \sim 2 (1+\sqrt{2}) \pi \max_{x \in \mathcal{W}(A)} \left\lvert \left(\frac{1 - x^{1/2^{s+1}}} {1 + x^{1/2^{s+1}}}  \right)^{2k+1} \right\rvert,
\]
for which we want to solve the problem of finding
\begin{equation}\label{eq:findingthevalues}
s,k \in \mathbb{N} \,:\; \|E(A,s,k)\| \sim 2 (1+\sqrt{2}) \pi \max_{x \in \mathcal{W}(A)} \left\lvert \left(\frac{1 - x^{1/2^{s+1}}} {1 + x^{1/2^{s+1}}}  \right)^{2k+1}  \right\rvert = \mathcal{E}(s,k) \leq \varepsilon,    
\end{equation}
for a given tolerance $\varepsilon$ on the final error; consider, e.g., the case depicted in Figure~\ref{fig:3dbound}. Since $s$ and $k$ are both integers, we can simply look at the different values of $\mathcal{E}(s,k)$. We also want to minimize the cost due to the computation of the square roots ($s \times \nicefrac{28}{3}n^3$ using \cite[Algorithm~6.7]{Higham}), and the evaluation of the Padé approximant ($k \times \nicefrac{2}{3}n^3$ using the LU factorization), therefore we can select
\[
(s,k) = \arg\min_{ (s,k)\,:\, \mathcal{E}(s,k) < \varepsilon } \left\lbrace \frac{28 s}{3} + \frac{2 k}{3} \right\rbrace = \arg\min_{ (s,k)\,:\, \mathcal{E}(s,k) < \varepsilon } c(s,k).
\]

A similar bound can be obtained also using the $\epsilon$-pseudospectra from Proposition~\ref{pro:pseudospectra} by substituting on the right hand-side the bound from Corollary~\ref{cor:bound-with-restriction}. {Succinctly, this allows us to choose $s$ and $k$ as 
\[
(s,k) = \arg\min_{ (s,k)\,:\, \mathcal{E}(s,k;\epsilon) < \varepsilon } \left\lbrace \frac{28 s}{3} + \frac{2 k}{3} \right\rbrace = \arg\min_{ (s,k)\,:\, \mathcal{E}(s,k;\epsilon) < \varepsilon } c(s,k), 
\]
for
\[
\mathcal{E}(s,k;\epsilon)   =   \frac{\mathcal{L}(\partial \Lambda_\epsilon(A))}{\epsilon}  \max_{x \in \Lambda_\epsilon(A)} \left\lvert \left( \frac{1- x^{1/2^{s+1}}   }{ 1+ x^{1/2^{s+1}}}\right)^{2k+1} \right\rvert.
\]
}
\section{Numerical experiments} \label{sec5}

The numerical examples are run with \texttt{MATLAB} version 9.12.0.1975300 (R2022a) on a Ubuntu 22.04.1 LTS (Jammy Jellyfish) notebook with 16 Gb of RAM and an Intel\textsuperscript{(R)} Core\textsuperscript{(TM)} i7-8750H CPU at 2.20 GHz. The codes to run the examples are available on the GitHub repository \href{https://github.com/Cirdans-Home/padelogarithm}{Cirdans-Home/padelogarithm}.

\subsection{Values of \texorpdfstring{$s$}{s} and \texorpdfstring{$k$}{k}}

Here we consider some example matrices and report the values of $s$ and $k$ chosen by the \texttt{MATLAB} \mintinline{matlab}{logm} routine~\cite{MR2970418,MR3080997} and through the bound~\eqref{eq:findingthevalues} together with the error of the approximation computed with respect to that of \mintinline{matlab}{logm}; see Table~\ref{tab:values_of_s_and_k}.

\begin{table}[htbp]
    \centering
    \caption{Number of inverse scaling and squaring steps $s$ made for different test matrices selected by the \mintinline{matlab}{logm} routine, and by using the bound~\eqref{eq:findingthevalues}. The first column reports the code needed to generate the test matrix. The last two columns report the absolute and relative errors computed with respect to the \mintinline{matlab}{logm} routine.}
    \label{tab:values_of_s_and_k}
    \small
    \begin{tabular}{p{0.35\textwidth}*{8}{c}}
    \toprule
    \multirow{2}{*}{Matrix} & \multicolumn{3}{c}{\mintinline{matlab}{logm}} & \multicolumn{3}{c}{Bound~\eqref{eq:findingthevalues}} &  \multirow{2}{*}{Abs. Err.} &  \multirow{2}{*}{Rel. Err.} \\
    & $s$ & $k$ & $c(s,k)$ & $s$ & $k$ & $c(s,k)$  \\
    \midrule
    \begin{minipage}[c]{0.4\textwidth}
\begin{minted}[fontsize=\footnotesize]{matlab}
A = gallery('forsythe',10,..
  1e-10,0)};
A = expm(A);
\end{minted}
\end{minipage} & 3  & 6 & 32  & 0 & 13 & \textbf{8} & 4.73e-16  & 4.73e-16 \\[2em]
    \begin{minipage}[c]{0.4\textwidth}
\begin{minted}[fontsize=\footnotesize]{matlab}
A = [cos(100) sin(100);
  -sin(100) cos(100)]
\end{minted}
\end{minipage} & 4 & 6 & 41 & 2 & 12 & \textbf{26} & 9.17e-16 & 3.24e-16 \\[2em]
\begin{minipage}[c]{0.4\textwidth}
\begin{minted}[fontsize=\footnotesize]{matlab}
n = 100;
A = gallery('triw',n,1,1);
A = A - diag(diag(A)) ...
+ diag(-(n-1)/2:(n-1)/2);
A = expm(A)
\end{minted}
\end{minipage} & 9 & 5 & 87 & 6 & 12 & \textbf{64} & 1.44e-14 & 2.89e-16 \\
    \bottomrule
    \end{tabular}
\end{table}

\subsection{Double exponential form}

We consider here a comparison with the adaptive double exponential (ADE) formula from \cite[Algorithm~2]{MR4074499}. We used always $k = 10$ node to start the adaptive step. The tolerances to be reached and the truncation has been set to achieve comparable absolute and relative errors with our approach. An estimate of the cost of the algorithm is given by the number of steps $k \times n^3$ for the ADE formula, and $(\frac{28 s}{3}+\frac{2 k}{3})\times n^3$ for the Gauss--Legendre with $s$ inverse scaling and squaring steps, and $k$ quadrature nodes.

\begingroup
\setlength{\tabcolsep}{0.45em}
\begin{table}[htbp]
    \centering
    \caption{Comparisons between the adaptive double exponential formula from \cite[Algorithm~2]{MR4074499} and the Gauss--Legendre rule with inverse scaling and squaring, and values of $s$ and $k$ fixed as in \eqref{eq:findingthevalues}.}
    \small
    \begin{tabular}{l*{8}{c}}
    \toprule
    & \multicolumn{3}{c}{ADE} & \multicolumn{5}{c}{Bound~\eqref{eq:findingthevalues}} \\
         Matrix & $k$ & Abs. Err. & Rel. Err. & $s$ & $k$ & $c(s,k)$ & Abs. Err. & Rel. Err. \\
\midrule
\mintinline[fontsize=\footnotesize]{matlab}{ gallery('parter',10);} & 73 & 5.39e-15 & 3.07e-15 & 1 & 12 & 17 & 8.06e-16 & 4.58e-16 \\
\mintinline[fontsize=\footnotesize]{matlab}{-gallery('hanowa',10);} & 145 & 1.49e-15 & 6.99e-16 & 1 & 15 & 19 & 4.97e-16 & 2.33e-16 \\
\mintinline[fontsize=\footnotesize]{matlab}{ gallery('dorr',10,0.05);} & 145 & 1.74e-14 & 5.20e-15 & 2 & 12 & 26 & 1.58e-15 & 4.71e-16 \\
    \bottomrule
    \end{tabular}
    
    \label{tab:adecomparison}
\end{table}
\endgroup
The error is computed again against the reference solution given by \mintinline{matlab}{logm}. The results in Table~\ref{tab:adecomparison} show that knowing \emph{a priori} 
the bound on the error permits to reduce of the computational cost. 

\subsection{\texorpdfstring{$\varepsilon$}{epsilon}-pseudospectrum bound}

{Let us consider here an example from~\cite[Section~3]{MR2200942}, for $a,b,c \in \mathbb{R}$ we consider the following parameterized family of matrices with known logarithm
\[
T_{a,b,c} = \exp(a) \begin{bmatrix}
    1 & b & \nicefrac{b^2}{2}+c \\
    0 & 1 & b \\
    0 & 0 & 1
\end{bmatrix}, \qquad \log\left(T_{a,b,c}\right) = \begin{bmatrix}
    a & b & c \\
    0 & a & b \\
    0 & 0 & a \\
\end{bmatrix}.
\] 
If we take values of $a \leq 0.05$, $b = 10^{3}$, and $c=10^{-3}$, then the resulting matrices have a field of values that exceeds the analyticity region of the function in the bound~\eqref{eq:bound-on-the-complex}. As we have discussed at the bottom of Section~\ref{sec4}, we can use in that case an $\epsilon$-pseudospectrum, and thus Proposition~\ref{pro:pseudospectra} in conjunction with~\eqref{eq:bound-on-the-complex} to determine the parameters.
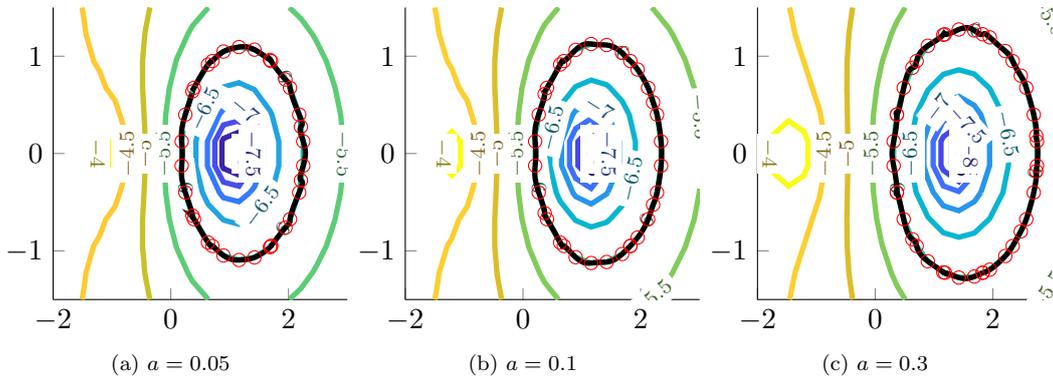
\begin{figure}[htb]
    \centering
    \subfloat[$a = 0.05$]{\begin{tikzpicture}

\begin{axis}[%
width=0.27\columnwidth,
height=0.27\columnwidth,
at={(1.236in,0.481in)},
scale only axis,
colormap={mymap}{[1pt] rgb(0pt)=(0.2422,0.1504,0.6603); rgb(1pt)=(0.2444,0.1534,0.6728); rgb(2pt)=(0.2464,0.1569,0.6847); rgb(3pt)=(0.2484,0.1607,0.6961); rgb(4pt)=(0.2503,0.1648,0.7071); rgb(5pt)=(0.2522,0.1689,0.7179); rgb(6pt)=(0.254,0.1732,0.7286); rgb(7pt)=(0.2558,0.1773,0.7393); rgb(8pt)=(0.2576,0.1814,0.7501); rgb(9pt)=(0.2594,0.1854,0.761); rgb(11pt)=(0.2628,0.1932,0.7828); rgb(12pt)=(0.2645,0.1972,0.7937); rgb(13pt)=(0.2661,0.2011,0.8043); rgb(14pt)=(0.2676,0.2052,0.8148); rgb(15pt)=(0.2691,0.2094,0.8249); rgb(16pt)=(0.2704,0.2138,0.8346); rgb(17pt)=(0.2717,0.2184,0.8439); rgb(18pt)=(0.2729,0.2231,0.8528); rgb(19pt)=(0.274,0.228,0.8612); rgb(20pt)=(0.2749,0.233,0.8692); rgb(21pt)=(0.2758,0.2382,0.8767); rgb(22pt)=(0.2766,0.2435,0.884); rgb(23pt)=(0.2774,0.2489,0.8908); rgb(24pt)=(0.2781,0.2543,0.8973); rgb(25pt)=(0.2788,0.2598,0.9035); rgb(26pt)=(0.2794,0.2653,0.9094); rgb(27pt)=(0.2798,0.2708,0.915); rgb(28pt)=(0.2802,0.2764,0.9204); rgb(29pt)=(0.2806,0.2819,0.9255); rgb(30pt)=(0.2809,0.2875,0.9305); rgb(31pt)=(0.2811,0.293,0.9352); rgb(32pt)=(0.2813,0.2985,0.9397); rgb(33pt)=(0.2814,0.304,0.9441); rgb(34pt)=(0.2814,0.3095,0.9483); rgb(35pt)=(0.2813,0.315,0.9524); rgb(36pt)=(0.2811,0.3204,0.9563); rgb(37pt)=(0.2809,0.3259,0.96); rgb(38pt)=(0.2807,0.3313,0.9636); rgb(39pt)=(0.2803,0.3367,0.967); rgb(40pt)=(0.2798,0.3421,0.9702); rgb(41pt)=(0.2791,0.3475,0.9733); rgb(42pt)=(0.2784,0.3529,0.9763); rgb(43pt)=(0.2776,0.3583,0.9791); rgb(44pt)=(0.2766,0.3638,0.9817); rgb(45pt)=(0.2754,0.3693,0.984); rgb(46pt)=(0.2741,0.3748,0.9862); rgb(47pt)=(0.2726,0.3804,0.9881); rgb(48pt)=(0.271,0.386,0.9898); rgb(49pt)=(0.2691,0.3916,0.9912); rgb(50pt)=(0.267,0.3973,0.9924); rgb(51pt)=(0.2647,0.403,0.9935); rgb(52pt)=(0.2621,0.4088,0.9946); rgb(53pt)=(0.2591,0.4145,0.9955); rgb(54pt)=(0.2556,0.4203,0.9965); rgb(55pt)=(0.2517,0.4261,0.9974); rgb(56pt)=(0.2473,0.4319,0.9983); rgb(57pt)=(0.2424,0.4378,0.9991); rgb(58pt)=(0.2369,0.4437,0.9996); rgb(59pt)=(0.2311,0.4497,0.9995); rgb(60pt)=(0.225,0.4559,0.9985); rgb(61pt)=(0.2189,0.462,0.9968); rgb(62pt)=(0.2128,0.4682,0.9948); rgb(63pt)=(0.2066,0.4743,0.9926); rgb(64pt)=(0.2006,0.4803,0.9906); rgb(65pt)=(0.195,0.4861,0.9887); rgb(66pt)=(0.1903,0.4919,0.9867); rgb(67pt)=(0.1869,0.4975,0.9844); rgb(68pt)=(0.1847,0.503,0.9819); rgb(69pt)=(0.1831,0.5084,0.9793); rgb(70pt)=(0.1818,0.5138,0.9766); rgb(71pt)=(0.1806,0.5191,0.9738); rgb(72pt)=(0.1795,0.5244,0.9709); rgb(73pt)=(0.1785,0.5296,0.9677); rgb(74pt)=(0.1778,0.5349,0.9641); rgb(75pt)=(0.1773,0.5401,0.9602); rgb(76pt)=(0.1768,0.5452,0.956); rgb(77pt)=(0.1764,0.5504,0.9516); rgb(78pt)=(0.1755,0.5554,0.9473); rgb(79pt)=(0.174,0.5605,0.9432); rgb(80pt)=(0.1716,0.5655,0.9393); rgb(81pt)=(0.1686,0.5705,0.9357); rgb(82pt)=(0.1649,0.5755,0.9323); rgb(83pt)=(0.161,0.5805,0.9289); rgb(84pt)=(0.1573,0.5854,0.9254); rgb(85pt)=(0.154,0.5902,0.9218); rgb(86pt)=(0.1513,0.595,0.9182); rgb(87pt)=(0.1492,0.5997,0.9147); rgb(88pt)=(0.1475,0.6043,0.9113); rgb(89pt)=(0.1461,0.6089,0.908); rgb(90pt)=(0.1446,0.6135,0.905); rgb(91pt)=(0.1429,0.618,0.9022); rgb(92pt)=(0.1408,0.6226,0.8998); rgb(93pt)=(0.1383,0.6272,0.8975); rgb(94pt)=(0.1354,0.6317,0.8953); rgb(95pt)=(0.1321,0.6363,0.8932); rgb(96pt)=(0.1288,0.6408,0.891); rgb(97pt)=(0.1253,0.6453,0.8887); rgb(98pt)=(0.1219,0.6497,0.8862); rgb(99pt)=(0.1185,0.6541,0.8834); rgb(100pt)=(0.1152,0.6584,0.8804); rgb(101pt)=(0.1119,0.6627,0.877); rgb(102pt)=(0.1085,0.6669,0.8734); rgb(103pt)=(0.1048,0.671,0.8695); rgb(104pt)=(0.1009,0.675,0.8653); rgb(105pt)=(0.0964,0.6789,0.8609); rgb(106pt)=(0.0914,0.6828,0.8562); rgb(107pt)=(0.0855,0.6865,0.8513); rgb(108pt)=(0.0789,0.6902,0.8462); rgb(109pt)=(0.0713,0.6938,0.8409); rgb(110pt)=(0.0628,0.6972,0.8355); rgb(111pt)=(0.0535,0.7006,0.8299); rgb(112pt)=(0.0433,0.7039,0.8242); rgb(113pt)=(0.0328,0.7071,0.8183); rgb(114pt)=(0.0234,0.7103,0.8124); rgb(115pt)=(0.0155,0.7133,0.8064); rgb(116pt)=(0.0091,0.7163,0.8003); rgb(117pt)=(0.0046,0.7192,0.7941); rgb(118pt)=(0.0019,0.722,0.7878); rgb(119pt)=(0.0009,0.7248,0.7815); rgb(120pt)=(0.0018,0.7275,0.7752); rgb(121pt)=(0.0046,0.7301,0.7688); rgb(122pt)=(0.0094,0.7327,0.7623); rgb(123pt)=(0.0162,0.7352,0.7558); rgb(124pt)=(0.0253,0.7376,0.7492); rgb(125pt)=(0.0369,0.74,0.7426); rgb(126pt)=(0.0504,0.7423,0.7359); rgb(127pt)=(0.0638,0.7446,0.7292); rgb(128pt)=(0.077,0.7468,0.7224); rgb(129pt)=(0.0899,0.7489,0.7156); rgb(130pt)=(0.1023,0.751,0.7088); rgb(131pt)=(0.1141,0.7531,0.7019); rgb(132pt)=(0.1252,0.7552,0.695); rgb(133pt)=(0.1354,0.7572,0.6881); rgb(134pt)=(0.1448,0.7593,0.6812); rgb(135pt)=(0.1532,0.7614,0.6741); rgb(136pt)=(0.1609,0.7635,0.6671); rgb(137pt)=(0.1678,0.7656,0.6599); rgb(138pt)=(0.1741,0.7678,0.6527); rgb(139pt)=(0.1799,0.7699,0.6454); rgb(140pt)=(0.1853,0.7721,0.6379); rgb(141pt)=(0.1905,0.7743,0.6303); rgb(142pt)=(0.1954,0.7765,0.6225); rgb(143pt)=(0.2003,0.7787,0.6146); rgb(144pt)=(0.2061,0.7808,0.6065); rgb(145pt)=(0.2118,0.7828,0.5983); rgb(146pt)=(0.2178,0.7849,0.5899); rgb(147pt)=(0.2244,0.7869,0.5813); rgb(148pt)=(0.2318,0.7887,0.5725); rgb(149pt)=(0.2401,0.7905,0.5636); rgb(150pt)=(0.2491,0.7922,0.5546); rgb(151pt)=(0.2589,0.7937,0.5454); rgb(152pt)=(0.2695,0.7951,0.536); rgb(153pt)=(0.2809,0.7964,0.5266); rgb(154pt)=(0.2929,0.7975,0.517); rgb(155pt)=(0.3052,0.7985,0.5074); rgb(156pt)=(0.3176,0.7994,0.4975); rgb(157pt)=(0.3301,0.8002,0.4876); rgb(158pt)=(0.3424,0.8009,0.4774); rgb(159pt)=(0.3548,0.8016,0.4669); rgb(160pt)=(0.3671,0.8021,0.4563); rgb(161pt)=(0.3795,0.8026,0.4454); rgb(162pt)=(0.3921,0.8029,0.4344); rgb(163pt)=(0.405,0.8031,0.4233); rgb(164pt)=(0.4184,0.803,0.4122); rgb(165pt)=(0.4322,0.8028,0.4013); rgb(166pt)=(0.4463,0.8024,0.3904); rgb(167pt)=(0.4608,0.8018,0.3797); rgb(168pt)=(0.4753,0.8011,0.3691); rgb(169pt)=(0.4899,0.8002,0.3586); rgb(170pt)=(0.5044,0.7993,0.348); rgb(171pt)=(0.5187,0.7982,0.3374); rgb(172pt)=(0.5329,0.797,0.3267); rgb(173pt)=(0.547,0.7957,0.3159); rgb(175pt)=(0.5748,0.7929,0.2941); rgb(176pt)=(0.5886,0.7913,0.2833); rgb(177pt)=(0.6024,0.7896,0.2726); rgb(178pt)=(0.6161,0.7878,0.2622); rgb(179pt)=(0.6297,0.7859,0.2521); rgb(180pt)=(0.6433,0.7839,0.2423); rgb(181pt)=(0.6567,0.7818,0.2329); rgb(182pt)=(0.6701,0.7796,0.2239); rgb(183pt)=(0.6833,0.7773,0.2155); rgb(184pt)=(0.6963,0.775,0.2075); rgb(185pt)=(0.7091,0.7727,0.1998); rgb(186pt)=(0.7218,0.7703,0.1924); rgb(187pt)=(0.7344,0.7679,0.1852); rgb(188pt)=(0.7468,0.7654,0.1782); rgb(189pt)=(0.759,0.7629,0.1717); rgb(190pt)=(0.771,0.7604,0.1658); rgb(191pt)=(0.7829,0.7579,0.1608); rgb(192pt)=(0.7945,0.7554,0.157); rgb(193pt)=(0.806,0.7529,0.1546); rgb(194pt)=(0.8172,0.7505,0.1535); rgb(195pt)=(0.8281,0.7481,0.1536); rgb(196pt)=(0.8389,0.7457,0.1546); rgb(197pt)=(0.8495,0.7435,0.1564); rgb(198pt)=(0.86,0.7413,0.1587); rgb(199pt)=(0.8703,0.7392,0.1615); rgb(200pt)=(0.8804,0.7372,0.165); rgb(201pt)=(0.8903,0.7353,0.1695); rgb(202pt)=(0.9,0.7336,0.1749); rgb(203pt)=(0.9093,0.7321,0.1815); rgb(204pt)=(0.9184,0.7308,0.189); rgb(205pt)=(0.9272,0.7298,0.1973); rgb(206pt)=(0.9357,0.729,0.2061); rgb(207pt)=(0.944,0.7285,0.2151); rgb(208pt)=(0.9523,0.7284,0.2237); rgb(209pt)=(0.9606,0.7285,0.2312); rgb(210pt)=(0.9689,0.7292,0.2373); rgb(211pt)=(0.977,0.7304,0.2418); rgb(212pt)=(0.9842,0.733,0.2446); rgb(213pt)=(0.99,0.7365,0.2429); rgb(214pt)=(0.9946,0.7407,0.2394); rgb(215pt)=(0.9966,0.7458,0.2351); rgb(216pt)=(0.9971,0.7513,0.2309); rgb(217pt)=(0.9972,0.7569,0.2267); rgb(218pt)=(0.9971,0.7626,0.2224); rgb(219pt)=(0.9969,0.7683,0.2181); rgb(220pt)=(0.9966,0.774,0.2138); rgb(221pt)=(0.9962,0.7798,0.2095); rgb(222pt)=(0.9957,0.7856,0.2053); rgb(223pt)=(0.9949,0.7915,0.2012); rgb(224pt)=(0.9938,0.7974,0.1974); rgb(225pt)=(0.9923,0.8034,0.1939); rgb(226pt)=(0.9906,0.8095,0.1906); rgb(227pt)=(0.9885,0.8156,0.1875); rgb(228pt)=(0.9861,0.8218,0.1846); rgb(229pt)=(0.9835,0.828,0.1817); rgb(230pt)=(0.9807,0.8342,0.1787); rgb(231pt)=(0.9778,0.8404,0.1757); rgb(232pt)=(0.9748,0.8467,0.1726); rgb(233pt)=(0.972,0.8529,0.1695); rgb(234pt)=(0.9694,0.8591,0.1665); rgb(235pt)=(0.9671,0.8654,0.1636); rgb(236pt)=(0.9651,0.8716,0.1608); rgb(237pt)=(0.9634,0.8778,0.1582); rgb(238pt)=(0.9619,0.884,0.1557); rgb(239pt)=(0.9608,0.8902,0.1532); rgb(240pt)=(0.9601,0.8963,0.1507); rgb(241pt)=(0.9596,0.9023,0.148); rgb(242pt)=(0.9595,0.9084,0.145); rgb(243pt)=(0.9597,0.9143,0.1418); rgb(244pt)=(0.9601,0.9203,0.1382); rgb(245pt)=(0.9608,0.9262,0.1344); rgb(246pt)=(0.9618,0.932,0.1304); rgb(247pt)=(0.9629,0.9379,0.1261); rgb(248pt)=(0.9642,0.9437,0.1216); rgb(249pt)=(0.9657,0.9494,0.1168); rgb(250pt)=(0.9674,0.9552,0.1116); rgb(251pt)=(0.9692,0.9609,0.1061); rgb(252pt)=(0.9711,0.9667,0.1001); rgb(253pt)=(0.973,0.9724,0.0938); rgb(254pt)=(0.9749,0.9782,0.0872); rgb(255pt)=(0.9769,0.9839,0.0805)},
unbounded coords=jump,
xmin=-2,
xmax=3,
ymin=-1.5,
ymax=1.5,
axis background/.style={fill=white},
axis x line*=bottom,
axis y line*=left
]
\addplot[contour prepared, contour prepared format=matlab, line width=2.0pt] table[] {%
pseudoa005-1.tsv};
\addplot [color=black, line width=2.0pt, forget plot]
  table[]{pseudoa005-2.tsv};
\addplot [color=black, line width=2.0pt, forget plot]
  table[]{pseudoa005-3.tsv};
\addplot [color=red, only marks, mark size=2.5pt, mark=o, mark options={solid, red}, forget plot]
  table[]{pseudoa005-4.tsv};
\end{axis}
\end{tikzpicture}%
}
    \subfloat[$a = 0.1$]{\begin{tikzpicture}

\begin{axis}[%
width=0.27\columnwidth,
height=0.27\columnwidth,
at={(1.236in,0.481in)},
scale only axis,
colormap={mymap}{[1pt] rgb(0pt)=(0.2422,0.1504,0.6603); rgb(1pt)=(0.2444,0.1534,0.6728); rgb(2pt)=(0.2464,0.1569,0.6847); rgb(3pt)=(0.2484,0.1607,0.6961); rgb(4pt)=(0.2503,0.1648,0.7071); rgb(5pt)=(0.2522,0.1689,0.7179); rgb(6pt)=(0.254,0.1732,0.7286); rgb(7pt)=(0.2558,0.1773,0.7393); rgb(8pt)=(0.2576,0.1814,0.7501); rgb(9pt)=(0.2594,0.1854,0.761); rgb(11pt)=(0.2628,0.1932,0.7828); rgb(12pt)=(0.2645,0.1972,0.7937); rgb(13pt)=(0.2661,0.2011,0.8043); rgb(14pt)=(0.2676,0.2052,0.8148); rgb(15pt)=(0.2691,0.2094,0.8249); rgb(16pt)=(0.2704,0.2138,0.8346); rgb(17pt)=(0.2717,0.2184,0.8439); rgb(18pt)=(0.2729,0.2231,0.8528); rgb(19pt)=(0.274,0.228,0.8612); rgb(20pt)=(0.2749,0.233,0.8692); rgb(21pt)=(0.2758,0.2382,0.8767); rgb(22pt)=(0.2766,0.2435,0.884); rgb(23pt)=(0.2774,0.2489,0.8908); rgb(24pt)=(0.2781,0.2543,0.8973); rgb(25pt)=(0.2788,0.2598,0.9035); rgb(26pt)=(0.2794,0.2653,0.9094); rgb(27pt)=(0.2798,0.2708,0.915); rgb(28pt)=(0.2802,0.2764,0.9204); rgb(29pt)=(0.2806,0.2819,0.9255); rgb(30pt)=(0.2809,0.2875,0.9305); rgb(31pt)=(0.2811,0.293,0.9352); rgb(32pt)=(0.2813,0.2985,0.9397); rgb(33pt)=(0.2814,0.304,0.9441); rgb(34pt)=(0.2814,0.3095,0.9483); rgb(35pt)=(0.2813,0.315,0.9524); rgb(36pt)=(0.2811,0.3204,0.9563); rgb(37pt)=(0.2809,0.3259,0.96); rgb(38pt)=(0.2807,0.3313,0.9636); rgb(39pt)=(0.2803,0.3367,0.967); rgb(40pt)=(0.2798,0.3421,0.9702); rgb(41pt)=(0.2791,0.3475,0.9733); rgb(42pt)=(0.2784,0.3529,0.9763); rgb(43pt)=(0.2776,0.3583,0.9791); rgb(44pt)=(0.2766,0.3638,0.9817); rgb(45pt)=(0.2754,0.3693,0.984); rgb(46pt)=(0.2741,0.3748,0.9862); rgb(47pt)=(0.2726,0.3804,0.9881); rgb(48pt)=(0.271,0.386,0.9898); rgb(49pt)=(0.2691,0.3916,0.9912); rgb(50pt)=(0.267,0.3973,0.9924); rgb(51pt)=(0.2647,0.403,0.9935); rgb(52pt)=(0.2621,0.4088,0.9946); rgb(53pt)=(0.2591,0.4145,0.9955); rgb(54pt)=(0.2556,0.4203,0.9965); rgb(55pt)=(0.2517,0.4261,0.9974); rgb(56pt)=(0.2473,0.4319,0.9983); rgb(57pt)=(0.2424,0.4378,0.9991); rgb(58pt)=(0.2369,0.4437,0.9996); rgb(59pt)=(0.2311,0.4497,0.9995); rgb(60pt)=(0.225,0.4559,0.9985); rgb(61pt)=(0.2189,0.462,0.9968); rgb(62pt)=(0.2128,0.4682,0.9948); rgb(63pt)=(0.2066,0.4743,0.9926); rgb(64pt)=(0.2006,0.4803,0.9906); rgb(65pt)=(0.195,0.4861,0.9887); rgb(66pt)=(0.1903,0.4919,0.9867); rgb(67pt)=(0.1869,0.4975,0.9844); rgb(68pt)=(0.1847,0.503,0.9819); rgb(69pt)=(0.1831,0.5084,0.9793); rgb(70pt)=(0.1818,0.5138,0.9766); rgb(71pt)=(0.1806,0.5191,0.9738); rgb(72pt)=(0.1795,0.5244,0.9709); rgb(73pt)=(0.1785,0.5296,0.9677); rgb(74pt)=(0.1778,0.5349,0.9641); rgb(75pt)=(0.1773,0.5401,0.9602); rgb(76pt)=(0.1768,0.5452,0.956); rgb(77pt)=(0.1764,0.5504,0.9516); rgb(78pt)=(0.1755,0.5554,0.9473); rgb(79pt)=(0.174,0.5605,0.9432); rgb(80pt)=(0.1716,0.5655,0.9393); rgb(81pt)=(0.1686,0.5705,0.9357); rgb(82pt)=(0.1649,0.5755,0.9323); rgb(83pt)=(0.161,0.5805,0.9289); rgb(84pt)=(0.1573,0.5854,0.9254); rgb(85pt)=(0.154,0.5902,0.9218); rgb(86pt)=(0.1513,0.595,0.9182); rgb(87pt)=(0.1492,0.5997,0.9147); rgb(88pt)=(0.1475,0.6043,0.9113); rgb(89pt)=(0.1461,0.6089,0.908); rgb(90pt)=(0.1446,0.6135,0.905); rgb(91pt)=(0.1429,0.618,0.9022); rgb(92pt)=(0.1408,0.6226,0.8998); rgb(93pt)=(0.1383,0.6272,0.8975); rgb(94pt)=(0.1354,0.6317,0.8953); rgb(95pt)=(0.1321,0.6363,0.8932); rgb(96pt)=(0.1288,0.6408,0.891); rgb(97pt)=(0.1253,0.6453,0.8887); rgb(98pt)=(0.1219,0.6497,0.8862); rgb(99pt)=(0.1185,0.6541,0.8834); rgb(100pt)=(0.1152,0.6584,0.8804); rgb(101pt)=(0.1119,0.6627,0.877); rgb(102pt)=(0.1085,0.6669,0.8734); rgb(103pt)=(0.1048,0.671,0.8695); rgb(104pt)=(0.1009,0.675,0.8653); rgb(105pt)=(0.0964,0.6789,0.8609); rgb(106pt)=(0.0914,0.6828,0.8562); rgb(107pt)=(0.0855,0.6865,0.8513); rgb(108pt)=(0.0789,0.6902,0.8462); rgb(109pt)=(0.0713,0.6938,0.8409); rgb(110pt)=(0.0628,0.6972,0.8355); rgb(111pt)=(0.0535,0.7006,0.8299); rgb(112pt)=(0.0433,0.7039,0.8242); rgb(113pt)=(0.0328,0.7071,0.8183); rgb(114pt)=(0.0234,0.7103,0.8124); rgb(115pt)=(0.0155,0.7133,0.8064); rgb(116pt)=(0.0091,0.7163,0.8003); rgb(117pt)=(0.0046,0.7192,0.7941); rgb(118pt)=(0.0019,0.722,0.7878); rgb(119pt)=(0.0009,0.7248,0.7815); rgb(120pt)=(0.0018,0.7275,0.7752); rgb(121pt)=(0.0046,0.7301,0.7688); rgb(122pt)=(0.0094,0.7327,0.7623); rgb(123pt)=(0.0162,0.7352,0.7558); rgb(124pt)=(0.0253,0.7376,0.7492); rgb(125pt)=(0.0369,0.74,0.7426); rgb(126pt)=(0.0504,0.7423,0.7359); rgb(127pt)=(0.0638,0.7446,0.7292); rgb(128pt)=(0.077,0.7468,0.7224); rgb(129pt)=(0.0899,0.7489,0.7156); rgb(130pt)=(0.1023,0.751,0.7088); rgb(131pt)=(0.1141,0.7531,0.7019); rgb(132pt)=(0.1252,0.7552,0.695); rgb(133pt)=(0.1354,0.7572,0.6881); rgb(134pt)=(0.1448,0.7593,0.6812); rgb(135pt)=(0.1532,0.7614,0.6741); rgb(136pt)=(0.1609,0.7635,0.6671); rgb(137pt)=(0.1678,0.7656,0.6599); rgb(138pt)=(0.1741,0.7678,0.6527); rgb(139pt)=(0.1799,0.7699,0.6454); rgb(140pt)=(0.1853,0.7721,0.6379); rgb(141pt)=(0.1905,0.7743,0.6303); rgb(142pt)=(0.1954,0.7765,0.6225); rgb(143pt)=(0.2003,0.7787,0.6146); rgb(144pt)=(0.2061,0.7808,0.6065); rgb(145pt)=(0.2118,0.7828,0.5983); rgb(146pt)=(0.2178,0.7849,0.5899); rgb(147pt)=(0.2244,0.7869,0.5813); rgb(148pt)=(0.2318,0.7887,0.5725); rgb(149pt)=(0.2401,0.7905,0.5636); rgb(150pt)=(0.2491,0.7922,0.5546); rgb(151pt)=(0.2589,0.7937,0.5454); rgb(152pt)=(0.2695,0.7951,0.536); rgb(153pt)=(0.2809,0.7964,0.5266); rgb(154pt)=(0.2929,0.7975,0.517); rgb(155pt)=(0.3052,0.7985,0.5074); rgb(156pt)=(0.3176,0.7994,0.4975); rgb(157pt)=(0.3301,0.8002,0.4876); rgb(158pt)=(0.3424,0.8009,0.4774); rgb(159pt)=(0.3548,0.8016,0.4669); rgb(160pt)=(0.3671,0.8021,0.4563); rgb(161pt)=(0.3795,0.8026,0.4454); rgb(162pt)=(0.3921,0.8029,0.4344); rgb(163pt)=(0.405,0.8031,0.4233); rgb(164pt)=(0.4184,0.803,0.4122); rgb(165pt)=(0.4322,0.8028,0.4013); rgb(166pt)=(0.4463,0.8024,0.3904); rgb(167pt)=(0.4608,0.8018,0.3797); rgb(168pt)=(0.4753,0.8011,0.3691); rgb(169pt)=(0.4899,0.8002,0.3586); rgb(170pt)=(0.5044,0.7993,0.348); rgb(171pt)=(0.5187,0.7982,0.3374); rgb(172pt)=(0.5329,0.797,0.3267); rgb(173pt)=(0.547,0.7957,0.3159); rgb(175pt)=(0.5748,0.7929,0.2941); rgb(176pt)=(0.5886,0.7913,0.2833); rgb(177pt)=(0.6024,0.7896,0.2726); rgb(178pt)=(0.6161,0.7878,0.2622); rgb(179pt)=(0.6297,0.7859,0.2521); rgb(180pt)=(0.6433,0.7839,0.2423); rgb(181pt)=(0.6567,0.7818,0.2329); rgb(182pt)=(0.6701,0.7796,0.2239); rgb(183pt)=(0.6833,0.7773,0.2155); rgb(184pt)=(0.6963,0.775,0.2075); rgb(185pt)=(0.7091,0.7727,0.1998); rgb(186pt)=(0.7218,0.7703,0.1924); rgb(187pt)=(0.7344,0.7679,0.1852); rgb(188pt)=(0.7468,0.7654,0.1782); rgb(189pt)=(0.759,0.7629,0.1717); rgb(190pt)=(0.771,0.7604,0.1658); rgb(191pt)=(0.7829,0.7579,0.1608); rgb(192pt)=(0.7945,0.7554,0.157); rgb(193pt)=(0.806,0.7529,0.1546); rgb(194pt)=(0.8172,0.7505,0.1535); rgb(195pt)=(0.8281,0.7481,0.1536); rgb(196pt)=(0.8389,0.7457,0.1546); rgb(197pt)=(0.8495,0.7435,0.1564); rgb(198pt)=(0.86,0.7413,0.1587); rgb(199pt)=(0.8703,0.7392,0.1615); rgb(200pt)=(0.8804,0.7372,0.165); rgb(201pt)=(0.8903,0.7353,0.1695); rgb(202pt)=(0.9,0.7336,0.1749); rgb(203pt)=(0.9093,0.7321,0.1815); rgb(204pt)=(0.9184,0.7308,0.189); rgb(205pt)=(0.9272,0.7298,0.1973); rgb(206pt)=(0.9357,0.729,0.2061); rgb(207pt)=(0.944,0.7285,0.2151); rgb(208pt)=(0.9523,0.7284,0.2237); rgb(209pt)=(0.9606,0.7285,0.2312); rgb(210pt)=(0.9689,0.7292,0.2373); rgb(211pt)=(0.977,0.7304,0.2418); rgb(212pt)=(0.9842,0.733,0.2446); rgb(213pt)=(0.99,0.7365,0.2429); rgb(214pt)=(0.9946,0.7407,0.2394); rgb(215pt)=(0.9966,0.7458,0.2351); rgb(216pt)=(0.9971,0.7513,0.2309); rgb(217pt)=(0.9972,0.7569,0.2267); rgb(218pt)=(0.9971,0.7626,0.2224); rgb(219pt)=(0.9969,0.7683,0.2181); rgb(220pt)=(0.9966,0.774,0.2138); rgb(221pt)=(0.9962,0.7798,0.2095); rgb(222pt)=(0.9957,0.7856,0.2053); rgb(223pt)=(0.9949,0.7915,0.2012); rgb(224pt)=(0.9938,0.7974,0.1974); rgb(225pt)=(0.9923,0.8034,0.1939); rgb(226pt)=(0.9906,0.8095,0.1906); rgb(227pt)=(0.9885,0.8156,0.1875); rgb(228pt)=(0.9861,0.8218,0.1846); rgb(229pt)=(0.9835,0.828,0.1817); rgb(230pt)=(0.9807,0.8342,0.1787); rgb(231pt)=(0.9778,0.8404,0.1757); rgb(232pt)=(0.9748,0.8467,0.1726); rgb(233pt)=(0.972,0.8529,0.1695); rgb(234pt)=(0.9694,0.8591,0.1665); rgb(235pt)=(0.9671,0.8654,0.1636); rgb(236pt)=(0.9651,0.8716,0.1608); rgb(237pt)=(0.9634,0.8778,0.1582); rgb(238pt)=(0.9619,0.884,0.1557); rgb(239pt)=(0.9608,0.8902,0.1532); rgb(240pt)=(0.9601,0.8963,0.1507); rgb(241pt)=(0.9596,0.9023,0.148); rgb(242pt)=(0.9595,0.9084,0.145); rgb(243pt)=(0.9597,0.9143,0.1418); rgb(244pt)=(0.9601,0.9203,0.1382); rgb(245pt)=(0.9608,0.9262,0.1344); rgb(246pt)=(0.9618,0.932,0.1304); rgb(247pt)=(0.9629,0.9379,0.1261); rgb(248pt)=(0.9642,0.9437,0.1216); rgb(249pt)=(0.9657,0.9494,0.1168); rgb(250pt)=(0.9674,0.9552,0.1116); rgb(251pt)=(0.9692,0.9609,0.1061); rgb(252pt)=(0.9711,0.9667,0.1001); rgb(253pt)=(0.973,0.9724,0.0938); rgb(254pt)=(0.9749,0.9782,0.0872); rgb(255pt)=(0.9769,0.9839,0.0805)},
unbounded coords=jump,
xmin=-2,
xmax=3,
ymin=-1.5,
ymax=1.5,
axis background/.style={fill=white},
axis x line*=bottom,
axis y line*=left
]
\addplot[contour prepared, contour prepared format=matlab, line width=2.0pt] table[] {%
pseudoa010-1.tsv};
\addplot [color=black, line width=2.0pt, forget plot]
  table[]{pseudoa010-2.tsv};
\addplot [color=black, line width=2.0pt, forget plot]
  table[]{pseudoa010-3.tsv};
\addplot [color=red, only marks, mark size=2.5pt, mark=o, mark options={solid, red}, forget plot]
  table[]{pseudoa010-4.tsv};
\end{axis}
\end{tikzpicture}%
}   
    \subfloat[$a = 0.3$]{\begin{tikzpicture}

\begin{axis}[%
width=0.27\columnwidth,
height=0.27\columnwidth,
at={(1.236in,0.481in)},
scale only axis,
colormap={mymap}{[1pt] rgb(0pt)=(0.2422,0.1504,0.6603); rgb(1pt)=(0.2444,0.1534,0.6728); rgb(2pt)=(0.2464,0.1569,0.6847); rgb(3pt)=(0.2484,0.1607,0.6961); rgb(4pt)=(0.2503,0.1648,0.7071); rgb(5pt)=(0.2522,0.1689,0.7179); rgb(6pt)=(0.254,0.1732,0.7286); rgb(7pt)=(0.2558,0.1773,0.7393); rgb(8pt)=(0.2576,0.1814,0.7501); rgb(9pt)=(0.2594,0.1854,0.761); rgb(11pt)=(0.2628,0.1932,0.7828); rgb(12pt)=(0.2645,0.1972,0.7937); rgb(13pt)=(0.2661,0.2011,0.8043); rgb(14pt)=(0.2676,0.2052,0.8148); rgb(15pt)=(0.2691,0.2094,0.8249); rgb(16pt)=(0.2704,0.2138,0.8346); rgb(17pt)=(0.2717,0.2184,0.8439); rgb(18pt)=(0.2729,0.2231,0.8528); rgb(19pt)=(0.274,0.228,0.8612); rgb(20pt)=(0.2749,0.233,0.8692); rgb(21pt)=(0.2758,0.2382,0.8767); rgb(22pt)=(0.2766,0.2435,0.884); rgb(23pt)=(0.2774,0.2489,0.8908); rgb(24pt)=(0.2781,0.2543,0.8973); rgb(25pt)=(0.2788,0.2598,0.9035); rgb(26pt)=(0.2794,0.2653,0.9094); rgb(27pt)=(0.2798,0.2708,0.915); rgb(28pt)=(0.2802,0.2764,0.9204); rgb(29pt)=(0.2806,0.2819,0.9255); rgb(30pt)=(0.2809,0.2875,0.9305); rgb(31pt)=(0.2811,0.293,0.9352); rgb(32pt)=(0.2813,0.2985,0.9397); rgb(33pt)=(0.2814,0.304,0.9441); rgb(34pt)=(0.2814,0.3095,0.9483); rgb(35pt)=(0.2813,0.315,0.9524); rgb(36pt)=(0.2811,0.3204,0.9563); rgb(37pt)=(0.2809,0.3259,0.96); rgb(38pt)=(0.2807,0.3313,0.9636); rgb(39pt)=(0.2803,0.3367,0.967); rgb(40pt)=(0.2798,0.3421,0.9702); rgb(41pt)=(0.2791,0.3475,0.9733); rgb(42pt)=(0.2784,0.3529,0.9763); rgb(43pt)=(0.2776,0.3583,0.9791); rgb(44pt)=(0.2766,0.3638,0.9817); rgb(45pt)=(0.2754,0.3693,0.984); rgb(46pt)=(0.2741,0.3748,0.9862); rgb(47pt)=(0.2726,0.3804,0.9881); rgb(48pt)=(0.271,0.386,0.9898); rgb(49pt)=(0.2691,0.3916,0.9912); rgb(50pt)=(0.267,0.3973,0.9924); rgb(51pt)=(0.2647,0.403,0.9935); rgb(52pt)=(0.2621,0.4088,0.9946); rgb(53pt)=(0.2591,0.4145,0.9955); rgb(54pt)=(0.2556,0.4203,0.9965); rgb(55pt)=(0.2517,0.4261,0.9974); rgb(56pt)=(0.2473,0.4319,0.9983); rgb(57pt)=(0.2424,0.4378,0.9991); rgb(58pt)=(0.2369,0.4437,0.9996); rgb(59pt)=(0.2311,0.4497,0.9995); rgb(60pt)=(0.225,0.4559,0.9985); rgb(61pt)=(0.2189,0.462,0.9968); rgb(62pt)=(0.2128,0.4682,0.9948); rgb(63pt)=(0.2066,0.4743,0.9926); rgb(64pt)=(0.2006,0.4803,0.9906); rgb(65pt)=(0.195,0.4861,0.9887); rgb(66pt)=(0.1903,0.4919,0.9867); rgb(67pt)=(0.1869,0.4975,0.9844); rgb(68pt)=(0.1847,0.503,0.9819); rgb(69pt)=(0.1831,0.5084,0.9793); rgb(70pt)=(0.1818,0.5138,0.9766); rgb(71pt)=(0.1806,0.5191,0.9738); rgb(72pt)=(0.1795,0.5244,0.9709); rgb(73pt)=(0.1785,0.5296,0.9677); rgb(74pt)=(0.1778,0.5349,0.9641); rgb(75pt)=(0.1773,0.5401,0.9602); rgb(76pt)=(0.1768,0.5452,0.956); rgb(77pt)=(0.1764,0.5504,0.9516); rgb(78pt)=(0.1755,0.5554,0.9473); rgb(79pt)=(0.174,0.5605,0.9432); rgb(80pt)=(0.1716,0.5655,0.9393); rgb(81pt)=(0.1686,0.5705,0.9357); rgb(82pt)=(0.1649,0.5755,0.9323); rgb(83pt)=(0.161,0.5805,0.9289); rgb(84pt)=(0.1573,0.5854,0.9254); rgb(85pt)=(0.154,0.5902,0.9218); rgb(86pt)=(0.1513,0.595,0.9182); rgb(87pt)=(0.1492,0.5997,0.9147); rgb(88pt)=(0.1475,0.6043,0.9113); rgb(89pt)=(0.1461,0.6089,0.908); rgb(90pt)=(0.1446,0.6135,0.905); rgb(91pt)=(0.1429,0.618,0.9022); rgb(92pt)=(0.1408,0.6226,0.8998); rgb(93pt)=(0.1383,0.6272,0.8975); rgb(94pt)=(0.1354,0.6317,0.8953); rgb(95pt)=(0.1321,0.6363,0.8932); rgb(96pt)=(0.1288,0.6408,0.891); rgb(97pt)=(0.1253,0.6453,0.8887); rgb(98pt)=(0.1219,0.6497,0.8862); rgb(99pt)=(0.1185,0.6541,0.8834); rgb(100pt)=(0.1152,0.6584,0.8804); rgb(101pt)=(0.1119,0.6627,0.877); rgb(102pt)=(0.1085,0.6669,0.8734); rgb(103pt)=(0.1048,0.671,0.8695); rgb(104pt)=(0.1009,0.675,0.8653); rgb(105pt)=(0.0964,0.6789,0.8609); rgb(106pt)=(0.0914,0.6828,0.8562); rgb(107pt)=(0.0855,0.6865,0.8513); rgb(108pt)=(0.0789,0.6902,0.8462); rgb(109pt)=(0.0713,0.6938,0.8409); rgb(110pt)=(0.0628,0.6972,0.8355); rgb(111pt)=(0.0535,0.7006,0.8299); rgb(112pt)=(0.0433,0.7039,0.8242); rgb(113pt)=(0.0328,0.7071,0.8183); rgb(114pt)=(0.0234,0.7103,0.8124); rgb(115pt)=(0.0155,0.7133,0.8064); rgb(116pt)=(0.0091,0.7163,0.8003); rgb(117pt)=(0.0046,0.7192,0.7941); rgb(118pt)=(0.0019,0.722,0.7878); rgb(119pt)=(0.0009,0.7248,0.7815); rgb(120pt)=(0.0018,0.7275,0.7752); rgb(121pt)=(0.0046,0.7301,0.7688); rgb(122pt)=(0.0094,0.7327,0.7623); rgb(123pt)=(0.0162,0.7352,0.7558); rgb(124pt)=(0.0253,0.7376,0.7492); rgb(125pt)=(0.0369,0.74,0.7426); rgb(126pt)=(0.0504,0.7423,0.7359); rgb(127pt)=(0.0638,0.7446,0.7292); rgb(128pt)=(0.077,0.7468,0.7224); rgb(129pt)=(0.0899,0.7489,0.7156); rgb(130pt)=(0.1023,0.751,0.7088); rgb(131pt)=(0.1141,0.7531,0.7019); rgb(132pt)=(0.1252,0.7552,0.695); rgb(133pt)=(0.1354,0.7572,0.6881); rgb(134pt)=(0.1448,0.7593,0.6812); rgb(135pt)=(0.1532,0.7614,0.6741); rgb(136pt)=(0.1609,0.7635,0.6671); rgb(137pt)=(0.1678,0.7656,0.6599); rgb(138pt)=(0.1741,0.7678,0.6527); rgb(139pt)=(0.1799,0.7699,0.6454); rgb(140pt)=(0.1853,0.7721,0.6379); rgb(141pt)=(0.1905,0.7743,0.6303); rgb(142pt)=(0.1954,0.7765,0.6225); rgb(143pt)=(0.2003,0.7787,0.6146); rgb(144pt)=(0.2061,0.7808,0.6065); rgb(145pt)=(0.2118,0.7828,0.5983); rgb(146pt)=(0.2178,0.7849,0.5899); rgb(147pt)=(0.2244,0.7869,0.5813); rgb(148pt)=(0.2318,0.7887,0.5725); rgb(149pt)=(0.2401,0.7905,0.5636); rgb(150pt)=(0.2491,0.7922,0.5546); rgb(151pt)=(0.2589,0.7937,0.5454); rgb(152pt)=(0.2695,0.7951,0.536); rgb(153pt)=(0.2809,0.7964,0.5266); rgb(154pt)=(0.2929,0.7975,0.517); rgb(155pt)=(0.3052,0.7985,0.5074); rgb(156pt)=(0.3176,0.7994,0.4975); rgb(157pt)=(0.3301,0.8002,0.4876); rgb(158pt)=(0.3424,0.8009,0.4774); rgb(159pt)=(0.3548,0.8016,0.4669); rgb(160pt)=(0.3671,0.8021,0.4563); rgb(161pt)=(0.3795,0.8026,0.4454); rgb(162pt)=(0.3921,0.8029,0.4344); rgb(163pt)=(0.405,0.8031,0.4233); rgb(164pt)=(0.4184,0.803,0.4122); rgb(165pt)=(0.4322,0.8028,0.4013); rgb(166pt)=(0.4463,0.8024,0.3904); rgb(167pt)=(0.4608,0.8018,0.3797); rgb(168pt)=(0.4753,0.8011,0.3691); rgb(169pt)=(0.4899,0.8002,0.3586); rgb(170pt)=(0.5044,0.7993,0.348); rgb(171pt)=(0.5187,0.7982,0.3374); rgb(172pt)=(0.5329,0.797,0.3267); rgb(173pt)=(0.547,0.7957,0.3159); rgb(175pt)=(0.5748,0.7929,0.2941); rgb(176pt)=(0.5886,0.7913,0.2833); rgb(177pt)=(0.6024,0.7896,0.2726); rgb(178pt)=(0.6161,0.7878,0.2622); rgb(179pt)=(0.6297,0.7859,0.2521); rgb(180pt)=(0.6433,0.7839,0.2423); rgb(181pt)=(0.6567,0.7818,0.2329); rgb(182pt)=(0.6701,0.7796,0.2239); rgb(183pt)=(0.6833,0.7773,0.2155); rgb(184pt)=(0.6963,0.775,0.2075); rgb(185pt)=(0.7091,0.7727,0.1998); rgb(186pt)=(0.7218,0.7703,0.1924); rgb(187pt)=(0.7344,0.7679,0.1852); rgb(188pt)=(0.7468,0.7654,0.1782); rgb(189pt)=(0.759,0.7629,0.1717); rgb(190pt)=(0.771,0.7604,0.1658); rgb(191pt)=(0.7829,0.7579,0.1608); rgb(192pt)=(0.7945,0.7554,0.157); rgb(193pt)=(0.806,0.7529,0.1546); rgb(194pt)=(0.8172,0.7505,0.1535); rgb(195pt)=(0.8281,0.7481,0.1536); rgb(196pt)=(0.8389,0.7457,0.1546); rgb(197pt)=(0.8495,0.7435,0.1564); rgb(198pt)=(0.86,0.7413,0.1587); rgb(199pt)=(0.8703,0.7392,0.1615); rgb(200pt)=(0.8804,0.7372,0.165); rgb(201pt)=(0.8903,0.7353,0.1695); rgb(202pt)=(0.9,0.7336,0.1749); rgb(203pt)=(0.9093,0.7321,0.1815); rgb(204pt)=(0.9184,0.7308,0.189); rgb(205pt)=(0.9272,0.7298,0.1973); rgb(206pt)=(0.9357,0.729,0.2061); rgb(207pt)=(0.944,0.7285,0.2151); rgb(208pt)=(0.9523,0.7284,0.2237); rgb(209pt)=(0.9606,0.7285,0.2312); rgb(210pt)=(0.9689,0.7292,0.2373); rgb(211pt)=(0.977,0.7304,0.2418); rgb(212pt)=(0.9842,0.733,0.2446); rgb(213pt)=(0.99,0.7365,0.2429); rgb(214pt)=(0.9946,0.7407,0.2394); rgb(215pt)=(0.9966,0.7458,0.2351); rgb(216pt)=(0.9971,0.7513,0.2309); rgb(217pt)=(0.9972,0.7569,0.2267); rgb(218pt)=(0.9971,0.7626,0.2224); rgb(219pt)=(0.9969,0.7683,0.2181); rgb(220pt)=(0.9966,0.774,0.2138); rgb(221pt)=(0.9962,0.7798,0.2095); rgb(222pt)=(0.9957,0.7856,0.2053); rgb(223pt)=(0.9949,0.7915,0.2012); rgb(224pt)=(0.9938,0.7974,0.1974); rgb(225pt)=(0.9923,0.8034,0.1939); rgb(226pt)=(0.9906,0.8095,0.1906); rgb(227pt)=(0.9885,0.8156,0.1875); rgb(228pt)=(0.9861,0.8218,0.1846); rgb(229pt)=(0.9835,0.828,0.1817); rgb(230pt)=(0.9807,0.8342,0.1787); rgb(231pt)=(0.9778,0.8404,0.1757); rgb(232pt)=(0.9748,0.8467,0.1726); rgb(233pt)=(0.972,0.8529,0.1695); rgb(234pt)=(0.9694,0.8591,0.1665); rgb(235pt)=(0.9671,0.8654,0.1636); rgb(236pt)=(0.9651,0.8716,0.1608); rgb(237pt)=(0.9634,0.8778,0.1582); rgb(238pt)=(0.9619,0.884,0.1557); rgb(239pt)=(0.9608,0.8902,0.1532); rgb(240pt)=(0.9601,0.8963,0.1507); rgb(241pt)=(0.9596,0.9023,0.148); rgb(242pt)=(0.9595,0.9084,0.145); rgb(243pt)=(0.9597,0.9143,0.1418); rgb(244pt)=(0.9601,0.9203,0.1382); rgb(245pt)=(0.9608,0.9262,0.1344); rgb(246pt)=(0.9618,0.932,0.1304); rgb(247pt)=(0.9629,0.9379,0.1261); rgb(248pt)=(0.9642,0.9437,0.1216); rgb(249pt)=(0.9657,0.9494,0.1168); rgb(250pt)=(0.9674,0.9552,0.1116); rgb(251pt)=(0.9692,0.9609,0.1061); rgb(252pt)=(0.9711,0.9667,0.1001); rgb(253pt)=(0.973,0.9724,0.0938); rgb(254pt)=(0.9749,0.9782,0.0872); rgb(255pt)=(0.9769,0.9839,0.0805)},
unbounded coords=jump,
xmin=-2,
xmax=3,
ymin=-1.5,
ymax=1.5,
axis background/.style={fill=white},
axis x line*=bottom,
axis y line*=left
]
\addplot[contour prepared, contour prepared format=matlab, line width=2.0pt] table[] {%
pseudoa030-1.tsv};
\addplot [color=black, line width=2.0pt, forget plot]
  table[]{pseudoa030-2.tsv};
\addplot [color=black, line width=2.0pt, forget plot]
  table[]{pseudoa030-3.tsv};
\addplot [color=red, only marks, mark size=2.5pt, mark=o, mark options={solid, red}, forget plot]
  table[]{pseudoa030-4.tsv};
\end{axis}
\end{tikzpicture}%
}
    \caption{$\epsilon$-pseudospectra for the matrices $T_{a,b,c}$ for different values of $a$. The thick black line ($\epsilon = 10^{-6}$) represents the level-set used for the evaluation of the bound with the red circles used as evaluation points. Since we only need it for the bound, the pseudospectrum is calculated with a grid of only twenty points.}
    \label{fig:cardoso-example}
\end{figure}
In Figure~\ref{fig:cardoso-example} we report the contour curves of the $\epsilon$-pseudospectra for different values of the parameter $a$ and $\epsilon = 10^{-6}$ which guarantee that the region in which we apply Proposition~\ref{pro:pseudospectra} is in the domain of analyticity of the error expression. To obtain such curves we used a combination of \texttt{EigTool}\footnote{See \href{https://github.com/eigtool/eigtool}{github.com/eigtool/eigtool}.}, for obtaining a \emph{grid sampling} of the resolvent function, and \texttt{Chebfun} to interpolate the level-set of interest and compute its arc-length.
\begin{table}[hbt]
    \centering
    \caption{{Bound and results obtained for the $\epsilon$-pseudospectra depicted in Figure~\ref{fig:cardoso-example}. We employ different values of $a$ and compare the absolute error in $\|\cdot\|_2$ between the exact logarithm $\log\left(T_{a,b,c}\right)$ and the approximated one obtained with the values of $s$ and $k$ from Proposition~\ref{pro:pseudospectra}. To compute the matrix logarithm the routines from the \texttt{advanpix} (\url{https://www.advanpix.com/}) library to work in quadruple precision have been used.}}
    \begin{tabular}{rllll}
    \toprule
       $a$  & $s$ & $k$ & Abs. Err. & Bound  \\
    \midrule
        0.05 & 2 & 11 & 1.03e-28 & 7.73e-17 \\
        0.1 & 2 & 11 & 1.59e-28 & 2.11e-17 \\
        0.3 & 2 & 10 & 3.07e-26 & 7.92e-18 \\
        0.5 & 2 & 10 & 5.04e-22 & 8.72e-18 \\
    \bottomrule
    \end{tabular}
    \label{tab:cardoso-example}
\end{table}
Indeed, we are using generic routines and better performance could be attained by using instead of a grid sampling algorithm a curve tracing one. In Table~\ref{tab:cardoso-example} we therefore report the different results in terms of error and the corresponding value of the bound. Observe that these test matrices are quite ill-conditioned, thus there is a limitation to the precision to which we can actually compute any matrix function on it. To show the properties of the bound, we therefore employ the \texttt{advanpix} library to perform the computation in augmented (quadruple) precision. The determination of the values of the pseudospectra, $s$ and $k$ was carried out with double precision as in all other cases. We observe that the bound is still predictive, but it is less accurate than in cases where we could use the field of value directly.}

\section{Conclusions} \label{sec6}

We have derived an accurate error analysis for a family of Padé approximations used for computing the logarithm of a matrix. This allows quantifying in a precise way the error committed with respect to the location of the spectrum of the matrices in question. The bound obtained can also be useful in describing the convergence of a rational Krylov method that uses as poles those of the rational approximation obtained through the Gauss--Legendre formulas. When information is available on the field of values of the matrix whose logarithm is to be evaluated, the new convergence analysis allows estimating \emph{a priori} the number of poles and inverse scaling and squaring steps necessary to achieve a certain accuracy.

\section*{Acknowledgement(s)}
We thank Massimiliano Fasi for a discussion on scaling and squaring algorithms. {We also thank the anonymous referee for pushing us to think further about the issues related to $\epsilon$-pseudospectra and the distance from normality of the matrices.}

\section*{Disclosure statement}

No potential conflict of interest was reported by the author(s).

\section*{Funding}
This work was partially supported by GNCS-INdAM. The authors are members of the INdAM research group GNCS.

\bibliographystyle{tfnlm}
\bibliography{bibliography}

\end{document}